\documentclass[final,reqno,a4paper]{siamltex}
\usepackage{stmaryrd,amssymb,amssymb,amsmath,enumerate}
\usepackage{epsfig,tabularx,color}

\newcommand{\eps}{\varepsilon}
\newcommand{\T}{\mathbb T}
\newcommand{\N}{\mathbb N}
\newcommand{\Z}{\mathbb Z}
\newcommand{\R}{\mathbb R}
\newcommand{\C}{\mathbb C}
\newcommand{\sch}{Schr\"{o}dinger }
\newcommand{\be}{\begin{equation}}
\newcommand{\ee}{\end{equation}}
\newcommand{\ba}{\begin{array}}
\newcommand{\ea}{\end{array}}

\newtheorem{remark}{Remark}[section]

\title{Solving highly-oscillatory NLS with SAM:\\ numerical efficiency and geometric properties
\thanks{The authors acknowledge support from the ANR-FWF Project LODIQUAS (ANR-11-IS01-0003).  
Norbert J. Mauser and Y. Zhang  acknowledge support from the Austrian Science Foundation (FWF) under grant No F41 (project VICOM), 
grant No I830 (project LODIQUAS) and the Austrian Ministry of Science
and Research via its grant for the WPI.}}

\author{Philippe Chartier\thanks{INRIA-Rennes, IRMAR and ENS Bruz, Campus de Beaulieu, 35042 Rennes Cedex, France, ({\tt Philippe.Chartier@inria.fr } (Ph. Chartier)).}
\and Norbert J. Mauser \thanks{Wolfgang Pauli Institute c/o Fak. Mathematik,
University Wien, Nordbergstrasse 15, 1090 Vienna, Austria ({\tt norbert.mauser@univie.ac.at
} (N.J. Mauser); {\tt yong.zhang@univie.ac.at} (Y. Zhang)).}
\and Florian M\'ehats \thanks{IRMAR, Universit\'e de Rennes 1 and INRIA-Rennes, Campus de Beaulieu, 
35042 Rennes Cedex, France ({\tt Florian.Mehats@univ-rennes1.fr} (F. M\'ehats)).}
\and Yong Zhang\footnotemark[3] }

\begin{document}

\maketitle

\begin{abstract} 
In this paper, we present the Stroboscopic Averaging Method (SAM), recently introduced in \cite{CCMSS11b,CCMSS11,SAM,CMSS10}, which aims at {\em numerically} solving highly-oscillatory differential equations. More specifically, we first apply SAM to the Schr\"odinger equation on the $1$-dimensional torus and on the real line with harmonic potential, with the aim of assessing its efficiency: as compared to the well-established standard splitting schemes, the stiffer the problem is, the larger the speed-up grows (up to a factor $100$ in our tests). The geometric properties of SAM are also explored: on very long time intervals, symmetric implementations of the method show a very good preservation of the mass invariant and of the energy. In a second series of experiments on $2$-dimensional equations, we demonstrate  the ability of SAM to capture qualitatively the long-time evolution of the solution (without spurring high oscillations).
\end{abstract}



\pagestyle{myheadings}\thispagestyle{plain}
 \markboth{Ph. Chartier, N. J. Mauser, F. M\'ehats and Y. Zhang}
{Solving highly-oscillatory NLS with SAM}

\section{Introduction}
This paper is devoted to the numerical solution of highly-oscillatory evolution equations posed in a Banach space
\begin{align} \label{eq:HOP}
\dot u^\eps (t)= \eps f\left(t,u^\eps(t)\right), \qquad t \in [0,T/\eps], \quad 
u^\eps(0) = u_0 \in X, 
\end{align}
where $f$ is a {\em smooth} map, periodic with respect to $t$, and where $T$ is a fixed positive time independent of $\eps$. The highly-oscillatory character of the equation stems from the length of the interval over which it is considered: this aspect is maybe better grasped whenever the problem is written in terms of the alternative rescaled variable $v^\eps(\tau) = u^\eps(t)$ with $\tau= \eps t$
$$
\frac{d}{d \tau} v^\eps(\tau) = \frac{1}{\eps} \frac{d}{d t} u^\eps(t) = f(\tau/\eps,v^\eps(\tau)), \qquad \tau \in [0,T], \quad 
v^\eps(0) = u_0, 
$$
a format for which $\eps$ clearly appears as the inverse of a frequency going to $\infty$ for $\eps$ going to $0$. 

In the two equivalent forms considered above, oscillations are {\em extrinsic}, i.e. explicitly present in the vector field $f$, as this is frequently the case in equations arising in practice and originating from physics or chemistry. Oscillations might also be {\em intrinsic}: in the application considered here, namely the Schr\"odinger equation, the problem in its original form is {\em autonomous} and no periodic variable is apparent. Nevertheless, it may be reformulated as in (\ref{eq:HOP}). Generally speaking, equations of the form
\begin{eqnarray} \label{eq:auton}
\dot v^\eps(t) = a(v^\eps(t)) + \eps b(v^\eps(t))
\end{eqnarray}
where the flow map $\chi_t$ of $\dot v = a(v)$ is periodic, can easily be rewritten, through the change of variable $v^\eps(t)=\chi_t(u^\eps(t))$, in  the form (\ref{eq:HOP}) with  
$$
f(t,u^\eps) = \Big(\frac{\partial \chi_t}{\partial u} (u^\eps) \Big)^{-1}b(\chi_t(u^\eps))
$$
and this transformation is extremely simple when $a$ is a linear map. \\

Problem (\ref{eq:HOP}) is notoriously difficult  to solve since numerical methods, in order to achieve some accuracy, are forced to follow more and more  oscillations (refer to (\ref{eq:HOP})) as $\eps$ becomes smaller and smaller, whereas the long-term dynamics is often what only matters in applications.  
Various numerical methods have been proposed in the literature for solving such problems. Standard methods such as Strang splitting or compositions thereof, which aim at solving (\ref{eq:auton}) directly, suffer from severe step size restrictions as $\eps$ goes to zero. More elaborated methods for (\ref{eq:auton}) when $a$ is  linear,  introduce filter functions in various ways (they are usually referred to as Gautschi-type methods) and bypass some of limitations of splitting techniques, although not in a completely satisfactory way as resonances are still present and/or geometric properties not preserved (see \cite{Hairer06}, Chapter XII for a survey).  The technique we present and use here proceeds in a completely different way. It relies upon the existence of an asymptotic high-order (in $\eps$) averaged equation associated with (\ref{eq:HOP}) and aims at approximating numerically the solution thereof through a micro-macro strategy. The underlying averaged equation 
\begin{eqnarray} \label{eq:Fi}
u^\eps(t) = \eps F^\eps(u^\eps) := \eps F_1(u^\eps) + \eps^2 F_2(u^\eps) + \ldots, \qquad u^\eps(0)=u_0,
\end{eqnarray}
being autonomous and in particular smooth w.r.t. $\eps$, it can be solved with a macro-integrator which benefits from the smallness of $\eps$, hence with a computational cost which is essentially independent of $\eps$ (the effect of the time-interval $T/\eps$ becoming larger with $\eps \rightarrow 0$ is indeed counterbalanced by the possibility to use also larger macro-steps $H={\cal O}(1/\eps)$). Evidently, $F^\eps$ needs to be computed in some way or another. Analytical expressions, though available, are becoming increasingly complex for high orders and are not easily amenable to practical computations. The strategy we use consists in  solving the original equation (\ref{eq:HOP}) with a micro-integrator over several periods and then combining the resulting values to approximate $F^\eps$ through finite difference formulas. Here, the stroboscopic character of SAM is crucial, as it allows to assert that the solutions of (\ref{eq:HOP}) and (\ref{eq:Fi}) coincide for values of $t$ that are multiple of the period. Altogether, SAM is a micro-macro procedure using only the original vector field $f$ and providing approximations of the exact solution of (\ref{eq:Fi}). We emphasize here again that its computational cost, in sharp contrast with standard integrators, is independent of $\eps$. A detailed description of stroboscopic averaging at both theoretical and practical levels shall be given in sections \ref{sect:sam} and \ref{sect:prop}.  \\

Another fundamental property of SAM is its geometric character. As a matter of  fact, it has been shown, first for ODEs \cite{CMSS10} and later on for evolution equations in a Hilbert space \cite{SAM}, that the averaged vector field $F^\eps$ in (\ref{eq:Fi}) inherits the properties of $f$. In particular, at a formal level, $F^\eps$ is Hamiltonian whenever $f$ is, and (\ref{eq:Fi}) has the same first integrals as $f$. For the case considered here, the quantities of interest are the Hamiltonian of the Schr\"odinger equation and the mass ($L^2$-norm) of the solution. These two quantities are still invariants of the associated averaged equation and it is of importance to design a version of SAM which also preserves these invariants numerically. This aspect will be further discussed in Section \ref{sect:prop}. \\

In this work, we apply SAM algorithm to the nonlinear Schr\"odinger equation 
$$
i\partial_t\psi^\eps(t,x)=(A\psi^\eps)(t,x)+\eps \alpha(x)|\psi^\eps(t,x)|^2\psi^\eps(t,x),\qquad t\geq 0 ,\quad
\psi(0,x)=\psi_0(x),
$$
either on the torus $\mathbb T_{2\pi}=[0,2\pi]$ with $A=-\partial^2_x$ and  periodic boundary conditions, or on $\R$  (Gross-Pitaevskii)  with $A=-\frac{1}{2}\partial^2_x+\frac{1}{2}(x^2-1)$. 
In both cases, the operator $A$ has its spectrum included in $\Z$ so that $e^{-itA}$ is $2\pi$-periodic in time. The equation can thus be regarded as an evolution equation of the form (\ref{eq:auton}) in a functional Hilbert space $X=H^s$. The first case is considered as a test equation to assess the efficiency of the method (see Subsection \ref{sect:torus}): we demonstrate  experimentally that the method can be as much as $100$ times faster than Strang splitting when $\eps$ is small. We then show in Subsection \ref{sect:torusgeo}, that whenever the macro-integrator is chosen appropriately, i.e. according to the discussion of Subsection \ref{sect:geonum}, SAM provides a numerical solution with constant mass along which the Hamiltonian is preserved up to a small error (which does not drift in time). The Gross-Pitaevskii equation in one dimension is used similarly to confirm the error analysis of  Subsection \ref{sect:prop}. \\

Finally, the objective  of Section \ref{sect:appli} is to 
illustrate how SAM can be used to explore the qualitative behavior of highly-oscillatory systems, in particular here the dynamics of the modes of the solution.  The nonlinearity in Schr\"odinger equations indeed induces coupling effects and energy transfers. To this aim, we shall consider three NLS models: Gross-Pitaevskii equation in one dimension and two models in two dimensions in the context of anisotropic confinement.

\section{Presentation of the stroboscopic averaging method} \label{sect:sam}

The so-called Stroboscopic Averaging Method (SAM in brief) was introduced in \cite{CCMSS11b,CCMSS11,CMSS10} for the purpose of solving highly-oscillatory ordinary differential equations. Its foundations rely on the asymptotic technique of stroboscopic averaging, whose aim is to write the exact solution of a differential system as the composition of a periodic (rapidly oscillating in time) change of variables with the flow of an autonomous (non-stiff) differential equation\footnote{The technique involves series that are generically not convergent, unless the problem is linear: in this particular case, it leads to a different though equivalent formulation of the well-known Floquet's theorem. In general, the sought decomposition can be obtained only up to a small remainder term.}. While various choices are possible for this change of variable, it is constructed in the framework of stroboscopic averaging so as to coincide with the identity map at  times that are multiple of the period, a property which will be crucial in the design of SAM. The relevance of this decomposition in the case of infinite dimensional Banach spaces was further analyzed  in \cite{SAM} and for application to the Schr\"odinger equation, we now present it in this context. 

\subsection{Stroboscopic change of variable and associated averaged vector field}
Given a differential equation, posed in a Banach space $X$, of the form
\begin{align} 
\dot u^\eps (t)= \eps f\left(t,u^\eps(t)\right), \qquad 
u^\eps(0) = u_0 \in X, 
\end{align}
where $f$ is a {\em smooth} map from $\T \times X$ into $X$ (with $\T \equiv \R/(P \Z)$) and where it it assumed that $u^\eps(t)$ exists on $[0,T/\eps]$ for some positive $T$, it can be shown that there exist a smooth change of variables
$\Phi^\eps: \T \times X \mapsto X$ and a smooth vector field $F^\eps: X \mapsto X$, such that 
\begin{align} \label{eq:decomp}
v^\eps(t) = \Phi^\eps \left(t,\Psi^\eps\left(t,u_0\right)\right)
\end{align} 
is an accurate approximation of the solution $u^\eps(t)$ of (\ref{eq:HOP}) in the sense that 
\begin{align} \label{eq:errorexp}
\forall\; t\in [0, T/\eps], \quad \|u^\eps(t) - v^\eps(t) \|_X \leq C e^{-C/\eps}
\end{align}
where $C$ is a positive constant independent of $\eps$. Here, $\Psi^\eps(t,u_0)$ denotes the $t$-flow of the differential equation 
\begin{align} \label{eq:aver}
\dot \Psi^\eps(t,u_0) = \eps F^\eps (\Psi^\eps(t,u_0)), \qquad 
\Psi^\eps(0,u_0) = u_0.
\end{align}
For such a result to hold, several assumptions are required, among which the most stringent one is the analyticity of $f$ with respect to its second argument $u \in X$. In contrast, only the continuity of $f$ w.r.t. its first argument is imposed. The next theorem (we refer to \cite{CMSS10} for a proof 
in the finite dimensional case and to \cite{SAM} in the context of Banach spaces) sums up the above discussion in precise mathematical terms: 
\begin{theorem} \label{th:main} Assume that {\bf (i)} $f$ is continuous w.r.t. its first variable, {\bf (ii)} $f$ is analytic w.r.t. to its second variable, {\bf (iii)} there exist  $\eps_0>0$ and a bounded open subset $K$ of $X$ such that $u^\eps(t)$ exists and remains in $K$ for all $0 < \eps < \eps_0$ and all $t \in [0,T/\eps]$ and {\bf (iv)} $f$ is bounded on $\T \times K$.

Then there exist  $\Phi^{\eps}$ (continuous and $P$-periodic w.r.t. its first variable, analytic w.r.t. to its second variable) and $F^\eps$  (analytic), and constants $0<\eps_1<\eps_0$ and  $C>0$, such that
\begin{align} \label{eq:factotoprout}
\forall \; 0 < \eps < \eps_1,\quad \forall\; t \in [0,T/\eps],
\quad
\left\|
u^\eps(t)
-
\Phi^{\eps}(t, \Psi^{\eps}(t, u_0))
\right\|_{X}
\leq
C \exp\left(-\frac{C}{\eps}\right),
\end{align}
where $\Psi^\eps$ is the flow of the differential equation
$$
\dot \Psi^\eps(t,u_0) = \eps F^\eps \left( \Psi^\eps(t,u_0)  \right).
$$
Furthermore, if assumptions {\bf (iii)} and {\bf (iv)} are valid with $T=+\infty$, then 
estimate (\ref{eq:factotoprout}) holds on $[0, \hat T/\eps^{1+\alpha}]$, for any $\hat T >0$ and any $0 < \alpha < 1$, with constants $\eps_1$ and $C$ now depending on $\alpha$ and $\hat T$.
\end{theorem}
\begin{remark}
It should be noticed that the dynamics of the weakly nonlinear equation (\ref{eq:HOP}) becomes non-trivial on intervals of length greater than $T/\eps$ for which the variation of $u^\eps(t)$ is expected to be of size  ${\cal O}(1)$. Estimate (\ref{eq:factotoprout}) itself is thus highly non-trivial on intervals of length $\hat T/\eps^{1+\alpha}$. 
\end{remark}
\begin{remark}
Whenever the function $f$ is not analytic but only of class $C^k$ w.r.t. its second variable, a weaker result still holds with an error estimate of the form 
\begin{align} \label{eq:errorpol}
\forall t\in [0, T/\eps], \quad \|u^\eps(t) - \Phi^{\eps}(t, \Psi^{\eps}(t, u_0)) \|_X \leq C \eps^k,
\end{align}
where $\Phi^\eps$ and $F^\eps$ are now  only $k$-times differentiable.
\end{remark}
\begin{remark}
\label{remark}
The so-called {\em averaged vector field} $F^\eps$ has an expansion in powers of $\eps$ of the form $F^\eps(u) = F_1(u) + \eps F_2(u) + \ldots$, whose terms are built up from $f$ and its derivatives: for instance, the first two terms read
$$
F_1(u) = \frac{1}{P} \int_0^P f(\tau,u) d\tau \quad \mbox{ and } \quad F_2(u) = -\frac{1}{2P} \int_0^P \int_0^\tau [f(s,u),f(\tau,u)] ds d\tau,
$$
where $[f,g]=(\partial_u f) g-(\partial_u g) f$ denotes the usual Lie-bracket of smooth functions.  The next terms in the expansion can also be explicitly written down using either the methodology used in \cite{SAM} or a formal nonlinear Floquet-Magnus expansion exposed in \cite{BCOR09}. Nevertheless, they require the resolution of complicated recurrence formulas.
\end{remark}
\subsection{Geometric aspects of stroboscopic averaging}
\label{sectgeom}
Assume now that $X$ is a Hilbert space (with scalar product $(\cdot,\cdot)_X$) and that it is densely  continuously embedded in some ambient Hilbert space $Z$ (with real scalar product $(\cdot,\cdot)_Z$). The vector field $f$ is said to be 
Hamiltonian if there exists a bounded map $J \in {\cal GL}(X)$, skew-symmetric w.r.t. $(\cdot,\cdot )_Z$, and an analytic  function $H: \T \times X \mapsto \C$  such that $f(t,u)=J^{-1}\nabla_u H(t,u)$  where the gradient is taken w.r.t. $(\cdot, \cdot)_Z$, i.e. 
$$\forall  (t,u,v),\qquad  (\nabla_u H(t,u) , v)_Z
=\partial_u H(t,u) \, v.$$
Accordingly, we say that $\Phi^\eps (t,u)$ is symplectic if
$$
\forall (t,u,v,w), \qquad \left( J \partial_u \Phi^\eps(t,u) v\,,\,\partial_u \Phi^\eps(t,u)w\right)_Z
= (J v,w)_Z.
$$
Finally, a smooth function $I^\eps: \T\times X \to \R$ (possibly depending on $\eps$) is said to be an invariant of (\ref{eq:HOP}) if it 
satisfies 
\begin{align}
\label{invava}
\forall (t,u), \quad \partial_t I^\eps(t,u) + \eps \, \partial_u I^\eps(t,u) f(t,u)=0,
\end{align}
which implies that 
$$
\forall t, \quad 
I^\eps(t,u^\eps(t))\equiv I^\eps(0,u_0).
$$
A number of important properties of stroboscopic averaging stem from  the fact that both the change of variable $\Phi^\eps$
and the averaged vector field $F^\eps$ inherit the intrinsic properties of the system. In particular (see \cite{SAM,CMSS10}):
\begin{itemize}
\item If the original equation (\ref{eq:HOP}) is Hamiltonian, then $\Phi^\eps$ is {\em symplectic} and $F^\eps$ is Hamiltonian. 
\item If (\ref{eq:HOP}) is divergence-free, then $\Phi^\eps$ is {\em volume-preserving} and $F^\eps$ is divergence-free. 
\item If $I^\eps: \T \times X \mapsto \C$ is an invariant of  (\ref{eq:HOP}) , then $I^\eps(0,\cdot): X \mapsto \C$ is preserved by $\Phi^\eps$ and is an invariant of $F^\eps$.
\end{itemize}
The second point is important in applications to kinetic equations for instance, although not relevant to our application to Schr\"odinger equation, for which it appears as a mere consequence of its Hamiltonian character. Hence, it will not be discussed further here.  \\

In order to later analyze the numerical experiments on a scientifically sound ground, we quote the following results from \cite{SAM}, also proved in \cite{CMSS10} in the finite-dimensional context:

\begin{theorem}
\label{prop:ham}
{\bf [Stroboscopic averaging preserves the Hamiltonian structure]}
Under the assumptions of Theorem \ref{th:main}, suppose that $f$ is Hamiltonian with Hamiltonian $H$. Then $\Phi^\eps$ and $F^\eps$ are respectively symplectic and  Hamiltonian up to exponentially small perturbation terms, namely there exists an analytic function $H^\eps: X\mapsto \C$ such that for all $(t,u,v,w) \in \T \times K \times X^2$:
\begin{eqnarray}
\label{sympl}
&\forall\; 0<\eps<\eps_1, \;
\Big|\left(
J \partial_u \Phi^\eps(t,u) v
\, , \,
\partial_u \Phi^\eps(t,u)w
\right)_Z-(J v \, , \, w )_Z \Big|
\leq C e^{-C/\eps}  \|v\|_X\,\|w\|_X,\quad \quad \\
 \label{ham} & \|F^\eps(u)-J^{-1} \nabla_u H^\eps(u)\|_X  \leq C e^{-C/\eps}.
\end{eqnarray}
\end{theorem}
\begin{remark} The Hamiltonian $H^\eps$ of the averaged equation, similarly to $F^\eps$, also has an expansion in $\eps$ of the form $H^\eps(u) = H_1(u) + \eps H_2(u) + \ldots$. Its terms are built up from $H$ and its derivatives: for instance, the first two terms read
$$
H_1(u) = \frac{1}{P} \int_0^P H(\tau,u) d\tau, \qquad H_2(u) = -\frac{1}{2P} \int_0^P \int_0^\tau\{H(s,u),H(\tau,u)\} ds d\tau
$$
where $\{\cdot,\cdot\}$ denotes the Poisson-bracket $\{H,G\}=(\nabla H, J \nabla G)_Z$ of smooth functions.\end{remark}
\begin{theorem}
\label{prop:inva}
{\bf [Stroboscopic averaging preserves invariants]}
Under the assumptions of Theorem \ref{th:main}, suppose
that the function $I^\eps: \T \times X \mapsto \R$ is an invariant of the field $f$, i.e. satisfies (\ref{invava}) for any $(t,u)\in\T\times K$ and is analytic w.r.t. its second variable. 
Then, $\Phi^\eps$ and $F^\eps$ satisfy for $0<\eps<\eps_1$ and $(t,u) \in \T \times K$
\begin{equation}
\quad \|I^\eps(t,\Phi^\eps(t,u)) - I^\eps(0,u) \|_X \leq C e^{-C/\eps}\;\mbox{and}\;\|\eps \partial_u I^\eps(0,u) \; F^\eps(u) \|_X \leq C e^{-C/\eps}
\end{equation}
for some positive constant $C$. 
In particular, we have
$$
\forall\; 0<\eps<\eps_1, \quad \forall \;t \in [0,T/\eps], \quad \|I^\eps(t,\Psi^\eps(t,u_0)) - I^\eps(0,u_0) \| \leq C e^{-C/\eps}.
$$
\end{theorem}
\begin{remark}
Both previous theorems still hold whenever the functions considered (i.e. $f$ and $I^\eps$) are only of class $C^k$. In this case, the $e^{-C/\eps}$-term has to replaced by $\eps^{k+1}$. 
\end{remark}

An important particular case arising in practice concerns the following   {\bf autonomous} semi-linear differential equation (a class to which Schr\"odinger equation belongs)
\begin{align} \label{eq:aut}
\dot w^\eps(t) = J^{-1} A w^\eps(t) + \eps g(w^\eps(t)), \qquad w^\eps(0)= w_0,
\end{align}
where $A$ is a linear (possibly unbounded) self-adjoint operator such that $e^{t J^{-1} A} \in {\cal L}(X,X)$ is P-periodic w.r.t. $t$. It is then straightforward to recast this system in the format of (\ref{eq:HOP}) by performing the preliminary change  
of variable (which is distinct from $\Phi^\eps$) $u^\eps(t)=e^{-t J^{-1} A} w^\eps(t)$, leading to 
\begin{equation}
\label{eqfiltree}
\dot u^\eps(t) = \eps \; e^{-t J^{-1} A} g\left( e^{t J^{-1} A} u^\eps(t) \right), \qquad u^\eps(0) = w^\eps(0).
\end{equation}
If additionally, if $g$ is derived from a potential $V$, i.e. if $g(w) = J^{-1} \nabla V(w)$, then equation (\ref{eq:aut}) is Hamitonian with Hamiltonian 
$$
\tilde H(w) = \frac{1}{2} (Aw, w)_Z+ \eps V(w),
$$
and so is equation (\ref{eq:HOP}) with $f(t,u) = e^{-t J^{-1} A} g\left( e^{t J^{-1} A} u \right)$ and 
$H(t,u) = V(e^{t J^{-1} A} u)$.
%
%
\subsection{SAM, a numerical counterpart of stroboscopic averaging}
The differential equation (\ref{eq:aver}) being non stiff, it makes good sense to try to approximate it rather than to solve the original stiff problem (\ref{eq:HOP}). However, a general-purpose numerical method can not rely on the analytical computation of such terms as those involved in $F^\eps$ and this rules out  the direct resolution of equation (\ref{eq:aver}). In the sequel, we rather solve it by approximating $F^\eps$ ``on the fly", which is the idea at the core of SAM and very much in the spirit of Heterogeneous Multiscale Methods (see \cite{E03,E03b,E05,E07}). In order to obtain an approximation of $F^\eps(u)$ at a given point $u \in X$, we first use the group property of $\Psi^\eps$ to assert that
$$
\eps F^\eps(u) = \left. \frac{d}{dt} \Psi^\eps(t,u)\right|_{t=0}
$$
and then approximate $F^\eps$ through an interpolation of the derivative of $\Psi^\eps(t,u)$, at order 2 by
\begin{align} \label{eq:inter}
F^\eps(u) \approx  \frac{1}{2P \eps} \left( \Psi^\eps(P,u) -  \Psi^\eps(-P,u) \right) = F_1(u) + \eps F_2(u) + {\cal O}(\eps^2)
\end{align}
or at order 4 by
\begin{align} \label{eq:inter2}
F^\eps(u) &\approx  \frac{1}{12P \eps} \left( -\Psi^\eps(2P,u)+8\Psi^\eps(P,u) - 8\Psi^\eps(-P,u)+\Psi^\eps(-2P,u) \right) \\
&= F_1(u) + \eps F_2(u)+\eps^2F_2(u) + \eps^3 F_3(u) + {\cal O}(\eps^4).\nonumber
\end{align}
To complete the procedure, it remains to use the fact that 
$$
 \Phi^\eps(P,\Psi^\eps(P,u)) = \Psi^\eps(P,u) \quad \mbox{ and } \quad \Phi^\eps(-P,\Psi^\eps(-P,u)) = \Psi^\eps(-P,u),
$$
a consequence of the stroboscopic property
$$
\forall \;k \in \Z, \quad \Phi^\eps(k P,\cdot) = {\rm Id}.
$$
It is worth insisting at this stage that the computation of $F^\eps$ at point $u$ is regarded as {\em asynchronous} in the terminology of Heterogeneous Multiscale Methods. This means here that $F^\eps(u)$ necessitates the computation of $\Psi^\eps(P,u)$ and 
$\Psi^\eps(-P,u)$ irrespectively of the time $t$ at which we compute the approximation of $v^\eps(t)$.  We finally obtain a numerical method by approximating $\Phi^\eps(P,\Psi^\eps(P,u))$ and $\Phi^\eps(-P,\Psi^\eps(-P,u))$ by solving 
the equations 
\begin{eqnarray}
 \label{eq:Fapprox}
&&\dot U^\eps = \eps f(t,U^\eps(u)), \; t \in [0,P], \;\;\; U^\eps(0)=u, \\[1ex]
&&\dot U^\eps = \eps f(t,U^\eps(u)), \; t \in [-P,0], \; U^\eps(0)=u,
\end{eqnarray}
by a standard one-step method $S^\eps_h$ (hereafter referred to as the {\em micro-integrator}) where the step size $h$ used is small enough to resolve one oscillation, i.e. $h=P/n$ with $n \in \N$. The outcome of this procedure is the following {\em micro-macro} algorithm: \\ \\
{\bf SAM Algorithm}
\begin{enumerate}
\item Choose a micro-step $h=P/n$ and a macro-step $H>0$ and set $N=0$.
\item Advance the solution through a standard explicit Runge-Kutta method (hereafter referred to as the {\em macro-integrator}) with coefficients $(a_{ij}, b_j)$:
\begin{align} \label{eq:RK}
u^i = u_{N} + \eps H \sum_{j=1}^{i-1} a_{ij} F^\eps_h(u^j),\; i=1 \ldots s, \;\; u_{N+1} = u_{N} + \eps H \sum_{j=1}^s b_j\, u^j,
\end{align}
where, either 
\begin{align} \label{eq:Fh}
F^\eps_h(u^j) = \frac{1}{2P\eps} \left((S^\eps_h)^n(u^j)-(S^\eps_{-h})^n(u^j) \right)
\end{align}
for second order interpolation, or
\begin{align} \label{eq:Fh2}
F^\eps_h(u^j) \!=\! \frac{1}{12P\eps} \left(-(S^\eps_h)^{2n}(u^j)+8(S^\eps_{h})^n(u^j)-8(S^\eps_{-h})^n(u^j)+(S^\eps_{-h})^{2n}(u^j) \right)
\end{align}
for fourth order interpolation.
\item Set $N:=N+1$ and go to step 2. until $N H \geq T/\eps$. 
\end{enumerate}

Note that the algorithm computes a sequence of approximations $u_N$ at times $t_N=NH$ to the averaged solution of (\ref{eq:aver}). For values of $t_N$ that coincide with integers multiple of the period $P$, this actually provides an approximation of $u^\eps(t_N)$, since $u^\eps(t_N)$ and $\Psi^\eps(t_N,u_0)$ then coincide. If one needs the solution at intermediate points, then it is also possible to obtain it through a kind of post-processing. If $kP < t_N < (k+1) P$ for some  $k \in \N$, an approximation of $u^\eps(t_N)$ is not directly  available, since $\Phi^\eps(t_N,\cdot) \neq {\rm Id}$. Nevertheless, it is straightforward to approximate $\Phi^\eps(t_N,u_N)$, since one has the relation:
\begin{align} \label{eq:inter3}
u^\eps(t_N) &= \Phi^\eps(t_N, \Psi^\eps(t_N,u_0))  = \Phi^\eps(\Delta t, \Psi^\eps(\Delta t, \Psi^\eps (-\Delta t,\Psi^\eps(t_N,u_0)))) \nonumber \\
&\approx \Phi^\eps(\Delta t, \Psi^\eps(\Delta t, \Psi^\eps (-\Delta t,u_N))) 
\end{align}
where we have used with $\Delta t = t_N-kP$ that, on the one hand 
$$
\Phi^\eps(t_N,\cdot) = \Phi^\eps(t_N-kP,\cdot) =  \Phi^\eps(\Delta t,\cdot) 
$$
owing to the periodicity of $\Psi^\eps$ w.r.t. $t$ and, on the other hand
$$
\Psi^\eps(t_N,\cdot) = \Psi^\eps(\Delta t, \Psi^\eps (-\Delta t,\Psi^\eps(t_N,\cdot)))
$$
owing to the group property of $\Psi^\eps$. Computing an approximation of $u^\eps(t_N)$ thus boils down to, first approximating $\tilde u^\eps_{kP} \approx \Psi^\eps(-\Delta t, u_N)$ by numerically  solving equation (\ref{eq:aver}) backward in time from $0$ to $-\Delta t$ and, second approximating $\Phi^\eps(\Delta t, \Psi^\eps(\Delta t, \tilde u^\eps_{kP}))$ by numerically solving equation (\ref{eq:HOP}) from $0$ to $\Delta t$ with initial value $u_0$ replaced by $\tilde u^\eps_{kP} $. These are only local-in-time operations which can be done independently of the main SAM algorithm and whose contribution to the global error of approximation is not significant.  In the following, we will not elaborate on this aspect of the method, since only stroboscopic times will be used in numerical simulations. 
%
%
\section{Numerical properties of SAM} \label{sect:prop}
\subsection{Formal error analysis}
\label{formal}
In this section, we present an error analysis of SAM. We shall content ourselves with a descriptive analysis, since a more mathematically rigorous derivation of error estimates  would unnecessarily burden the presentation, while bringing no surprise: standard arguments for non-stiff equations indeed hold, given that SAM requires to numerically solve the averaged equation (which is by construction non-stiff) and the original highly-oscillatory equation with step size $h$ much smaller than a period (thus getting the equation  back into a non-stiff regime). As explained in \cite{CCMSS11}, there are three sources of errors:
\begin{enumerate}
\item The approximation of $\eps F^\eps$ by a finite-difference formula (\ref{eq:inter}) or (\ref{eq:inter2}): the corresponding error is expressed in terms of $\eps$.
\item The substitution of $\Psi^\eps(\pm \ell P,u)$, $\ell=1,2$ in $F^\eps$ by $(S^\eps_{\pm h})^{\ell \cdot n}(u)$ in $F^\eps_h$ (\ref{eq:Fh}) or (\ref{eq:Fh2}): the corresponding error is expressed in terms of the micro-step $h$ and $\eps$.
\item The discretization error due to the macro-integrator (\ref{eq:RK}): the corresponding error is expressed in terms of the macro-step $H$.
\end{enumerate}
The first source of error contributes to ${\cal O}(\eps^{\delta+1})$ where $\delta$ is the order of the difference formula (either $2$ or $4$ in our experiments). The second source of error contributes to ${\cal O}(\eps^\nu h^p)$ since one solves a non-stiff equation with vector field $\eps f$ and step size $h$ over a bounded interval\footnote{The integer $p$ is here the order of the micro-integrator.}. The extra-factor $\eps^\nu$ for $\nu=1$ accounts for the fact that the micro-integrator is exact for $\eps=0$. It will turn out that $\nu=2$ is the effective value observed\footnote{The explanation of this behavior is rather nontrivial and somehow orthogonal to SAM. We refer to \cite{superconvergence} for a detailed study of this exponent $2$.}. These errors are magnified by a factor $1/\eps$ through the  stable macro-integration over an interval of length $T/\eps$ and lead to an error in $u_N$ of size ${\cal O}(\eps^{\delta}+\eps^{\nu-1} h^p)$.  The third source of error accounts for an error term of the form ${\cal O}((\eps H)^P)$ where $P$ is the order of the macro-integrator. We recall that $\eps H$ has the form $T/N$. \\ \\
Combined together, these three sources of error lead to an error of size 
\begin{equation}
\label{error}
{\cal O}(\eps^{\delta}+\eps^{\nu-1} h^p+(\eps H)^P).
\end{equation}
\begin{remark}
In the specific case of equation (\ref{eq:aut}), it is worth mentioning that the computation of $\Psi^\eps(kP,u)$ for $k=\pm 1, \pm 2$ for instance, can be done directly by solving the original equation (\ref{eq:aut}) (instead of  the filtered equation (\ref{eqfiltree})) on $1$ or $2$ periods, forward or backward in time. As a matter of fact, the filtering operation $u^\eps(t) = e^{-tJ^{-1} A} w_\eps(t)$ is transparent at times that are multiple of the period. 
\end{remark}
\subsection{Geometric behavior} 
\label{sect:geonum}
Although the asymptotic averaged field (\ref{eq:aver}) inherits geometric properties of the original system (\ref{eq:HOP}), there is no guarantee that its numerical implementation will also do so. When applied to a Hamiltonian equation such as the Schr\"odinger equation, the averaged vector field is also Hamiltonian and it may seem desirable for the numerical counterpart to be so. Unfortunately, even if  $(S_h)^n$ is a symplectic map, the finite-difference approximation is not a geometric transformation and our implementation of SAM is not symplectic. However, time-reversibility of the system (\ref{eq:HOP}), --whenever it holds--, is preserved provided the micro-integrator is itself a symmetric method (this will be the case of the splitting methods used in our experiments). Since in addition, $F^\eps_h$ is computed through symmetric formulas, then it is also the vector field of a time-reversible equation. This property ensures a favorable numerical behavior, as documented for ordinary-differential equations (see for instance \cite{Hairer06}): properties such as the preservation of the Hamitonian or the persistence of KAM-tori are indeed transferred through a symmetric discretization. It is thus in principle not difficult to design a symmetric SAM by choosing a symmetric macro-integrator. 

In the sequel, we use two kinds of macro-integrators. For simulations on time intervals of order $1/\eps$, our aim is to keep SAM as efficient as possible for small values of the parameter $\eps$ (i.e. for very high-oscillatory) for which SAM is expected to outperform existing methods, so we use non symmetric but accurate macro-integrators such as fourth-order Runge-Kutta method (RK4). On such intervals of time, the non stiff character of equation (\ref{eq:aver}) is enough for the numerical solution to exhibit essentially no drift in the energy and mass. On longer times though (such as $1/\eps^\beta$ with $\beta \geq 2$), we use a (second-order) symmetric method for the macro-integration, which exhibits much better preservation of geometric invariants.
%
%

\section{Numerical assessment tests}
\label{sect:num1}
The goal of this section is to confirm by numerical tests the error analysis sketched in Subsection \ref{formal} and to observe in addition the geometric properties of SAM for long-time integration. We apply SAM to nonlinear Schr\"odinger equations when the linear evolution induced by the Laplacian is periodic in time, see \cite{SAM}. We also compare the efficiency of SAM to the one of the time-splitting method. 

We consider cubic nonlinear Schr\"odinger equations of the form
\begin{align}
\label{NLS}
i\partial_t\psi^\eps(t,x)&=(A\psi^\eps)(t,x)+\eps \;\alpha(x)|\psi^\eps(t,x)|^2\psi^\eps(t,x),\qquad t\geq 0,\quad x\in \Omega,\\
\psi(0,x)&=\psi_0(x),\nonumber
\end{align}
in the two following situations:
\begin{enumerate}[(i)]
\item (NLS on the torus) $\Omega=\mathbb T_{2\pi}=[0,2\pi]$, $A=-\partial^2_x$ with periodic boundary conditions and $\alpha(x)=2\cos (2x)$; this problem was studied in \cite{Grebert}. We will take initial data which belong to the domain of $A$:
$$X=D(A)=\left\{u\in H^2(\mathbb T_{2\pi})\,:\, u(2\pi)=u(0)\right\}.$$
\item (Gross-Pitaevskii) $\Omega=\R$, $A=-\frac{1}{2}\partial^2_x+\frac{1}{2}(x^2-1)$ and $\alpha(x)\equiv 1$. Here also, we will choose initial data in the domain of $A$:
$$X=D(A)=\left\{u\in H^2(\R)\,:\, x^2\,u\in L^2(\R)\right\}.$$
\end{enumerate}
In both cases, the pivot space introduced in Subsection \ref{sectgeom} is $Z=L^2(\Omega)$. Moreover, the operator $A$ has compact resolvent and its spectrum is $\Z$ (NLS on the torus) or $\N$ (Gross-Pitaevskii). Hence, \eqref{NLS} is under the form of the autonomous equation \eqref{eq:aut} and the operator $e^{-itA}\in\mathcal L(X,X)$ is $2\pi$-periodic in time.

\subsection{NLS on the torus: accuracy and efficiency of SAM}\label{sect:torus}
We present here our numerical results in the case (i):
\begin{align}\label{micro:periodic}
i\partial_t\psi^\varepsilon&=  -\partial_{xx}\psi^\varepsilon + 2\varepsilon  \cos(2x)  |\psi^\varepsilon|^2 \psi^\varepsilon, \qquad 0 \le t \le T_0/ \varepsilon,\quad x\in \mathbb T_{2\pi},\\
\label{micro:periodic:ini}
\psi^\varepsilon(0,x) &= \psi_0(x)=\cos(x)+\sin(x),  \qquad x\in {\mathbb T_{2\pi}}.
\end{align}
As micro-integrator, we adopt a {\sl time-splitting spectral method}, which has been widely and successfully used in many applications \cite{BaoJM,BaoJM2}. More precisely, we choose the fourth-order time-splitting Fourier spectral method (TSFP4) introduced in \cite{Yoshida}. In practice, the wavefunction is discretized in space by Fourier series as
$\psi(x) = \sum_{-N_x/2+1}^{N_x/2} \widehat{\psi}_k e^{ikx}$, with $N_x=256$, such that errors originating from space discretization can be considered as negligible. The time-step for the micro-integrator is denoted by $h$ and is always taken under the form $2\pi/N$, where $N\in \N^*$. This ensures super-convergence, {\em i.e.} \eqref{error} holds with the exponent $\nu=2$ in , as proven in \cite{superconvergence}.

As macro-integrator, we use RK4 scheme with time-step $H=T_0/(\widetilde N \eps)$ and we approximate the vector field $F^\eps_h$ by using the fourth-order interpolation formula \eqref{eq:Fh2}. From \eqref{error}, we thus expect an asymptotic error estimate  of the form
\begin{equation}
\label{error2}
\mbox{\em error}_{\mbox{\scriptsize SAM}}\approx  C_1\eps^4+C_2\eps h^4+C_3(\eps H)^4.
\end{equation}
The final time of the simulation is taken as $T_0/\eps$ with $T_0=\pi/4$, while numerical errors are defined by 
$$\mbox{\em error} = \|\psi^{ref}-\psi^n\|_{\ell^2}=\left(\sum_{j=0}^{N_x-1}|\psi^{ref}_j -\psi_j^n|^2 \Delta x\right)^{1/2},$$ 
where $\psi^n$ denotes the numerical  solution. Finally, we take the value $\Delta x = 2\pi/N_x$ for the mesh size. For $\varepsilon = 2^{-3},\ldots,2^{-8}$, the reference solution $\psi^{ref}$ is obtained by the TSFP4 method on the whole time interval $[0,T_0/\eps]$, with a time-step $h=\eps \pi/2^{14}$. For  $\varepsilon = 2^{-9},2^{-10}$ (and smaller $\varepsilon$), the reference solution is obtained by SAM with $(\eps H, h) =(\pi/2^{12}, \pi/2^{12})$ using a higher-order interpolation method for the vector field (8th-order interpolation).
 
\begin{figure}[!htbp]
\centerline{
\psfig{figure=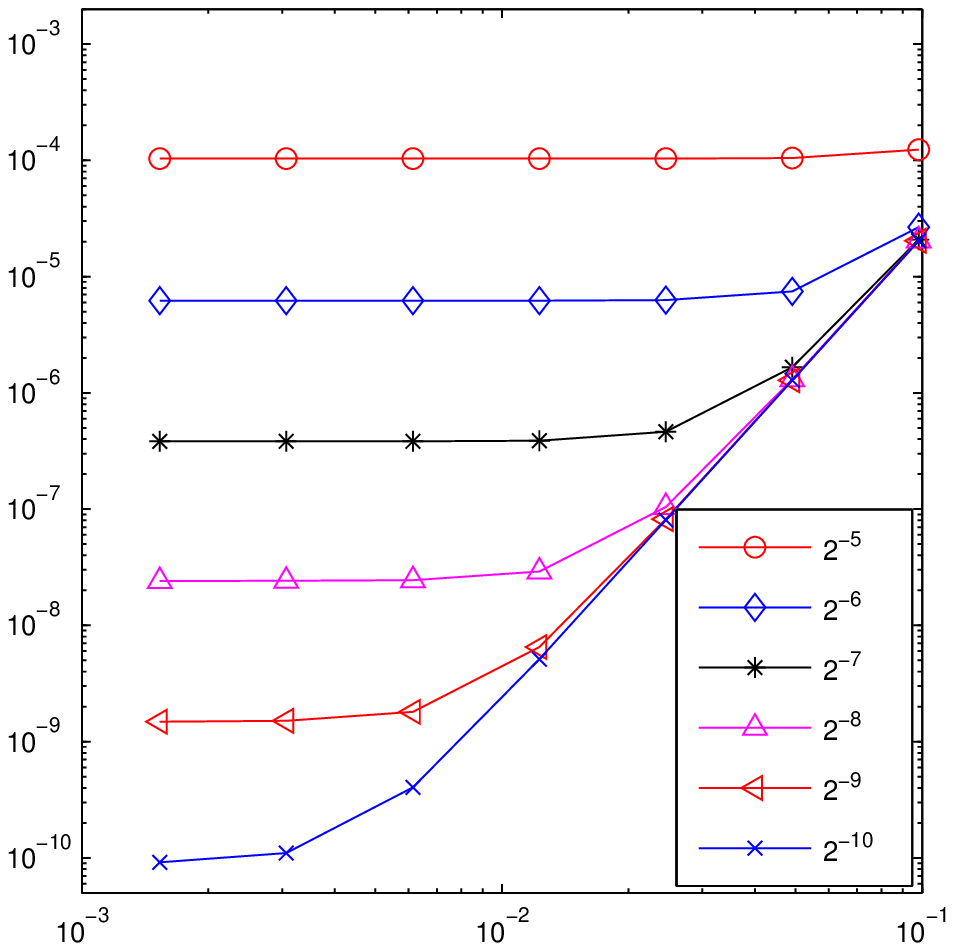,height=7cm,width=7.5cm}
\psfig{figure=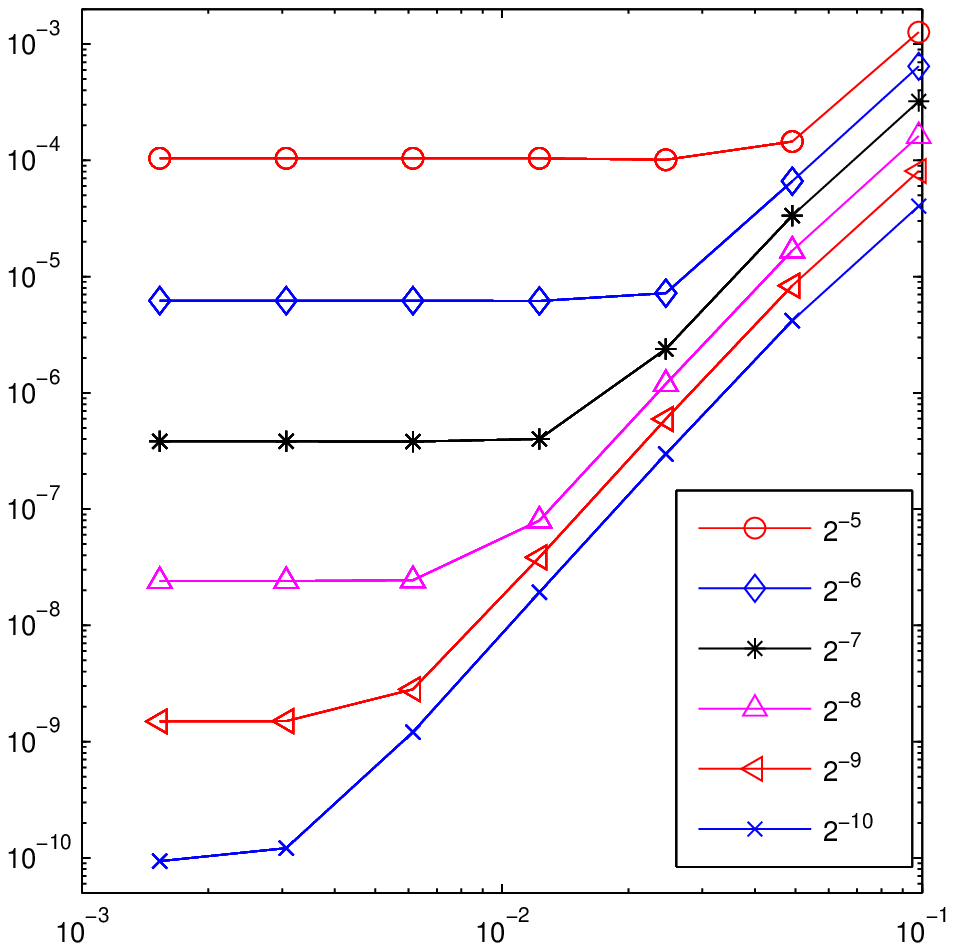,height=7cm,width=7.5cm}}
\caption{(NLS on the torus) Accuracy diagrams of SAM. Error versus macro step $H$ for $\eps=2^{-j/2}$, $j\in\{5,6,7,8,9,10\}$ (left) and error versus micro step $h$ for the same values of $\eps$ (right).} 
\label{fig:rk4_tsfp4}
\end{figure}

Our first results concern the accuracy of SAM and confirm estimate \eqref{error2}. On the left of  Fig. \ref{fig:rk4_tsfp4}, we present the error versus the macro step for  $H=\pi/(2^{j}\eps)$, $j = 5,6,\ldots,11$, with small fixed micro step $h =\pi/2^{12}$ The different curves correspond to the values $\eps=2^{-5},\ldots,2^{-10}$. On the right of  Fig. \ref{fig:rk4_tsfp4}, we represent the error versus the micro-step, for $h =\pi/2^{5},\ldots,\pi/2^{11}$, with small fixed macro step $H=\pi/(2^{12}\eps)$. Again, the different curves correspond to the values $\eps=2^{-5},\ldots,2^{-10}$.  Fourth-order accuracy in $\eps H$ (with a uniform constant w.r.t. $\eps$) and $h$ (with a linear constant in $\eps$) can be clearly observed. On each curve, a saturation appears, due to the interpolation error. As expected, the level of this saturation error is proportional to $\eps^4$.

Next, we investigate numerically the accuracy of the TSFP4 method. Table \ref{TSFP4_Accuracy} lists errors obtained by TSFP4 with time-step $h$ for different $\varepsilon$. We obtain an error of the form 
\begin{equation}
\label{errTSFP4}
\mbox{error}_{\mbox{\scriptsize TSFP4}}\approx C_4\varepsilon h^4,
\end{equation}
as expected from \cite{superconvergence}, since the micro time step is a submultiple of the period $2\pi$.
\begin{table}[t!]
\tabcolsep 0pt \caption{Errors of TSFP4 for NLS on $\mathbb T_{2\pi}$}
\label{TSFP4_Accuracy}\vspace*{-1em}
\begin{center}
\def\temptablewidth{1\textwidth}
{\rule{\temptablewidth}{1pt}}
\begin{tabular*}{\temptablewidth}{@{\extracolsep{\fill}}llllllll}
&$\varepsilon = 2^{-3}$ & $\varepsilon = 2^{-4}$&$\varepsilon = 2^{-5}$&$\varepsilon= 2^{-6}$ &$\varepsilon= 2^{-7}$&$\varepsilon = 2^{-8}$&$\varepsilon = 2^{-9}$\\
$h=\pi/2^4$ &6.17E-02&3.09E-02&1.54E-02&7.71E-03&3.86E-03&1.93E-03&9.64E-04\\
$h=\pi/2^5$ &5.50E-03&2.63E-03&1.30E-03&6.47E-04&3.23E-04&1.62E-04&8.08E-05\\
$h=\pi/2^6$ &6.15E-04&2.78E-04&1.35E-04&6.71E-05&3.35E-05&1.67E-05&8.37E-06\\
$h=\pi/2^7$ &4.31E-05&1.97E-05&9.61E-06&4.78E-06&2.38E-06&1.19E-06&5.96E-07\\
$h=\pi/2^8$ &2.86E-06&1.28E-06&6.22E-07&3.09E-07&1.54E-07&7.70E-08&3.85E-08\\
$h=\pi/2^9$ &1.82E-07&8.10E-08&3.92E-08&1.95E-08&9.71E-09&4.88E-09&2.43E-09\\
\end{tabular*}
{\rule{\temptablewidth}{1pt}}
\end{center}
\end{table}

\begin{figure}[!htbp]
\centerline{\psfig{figure=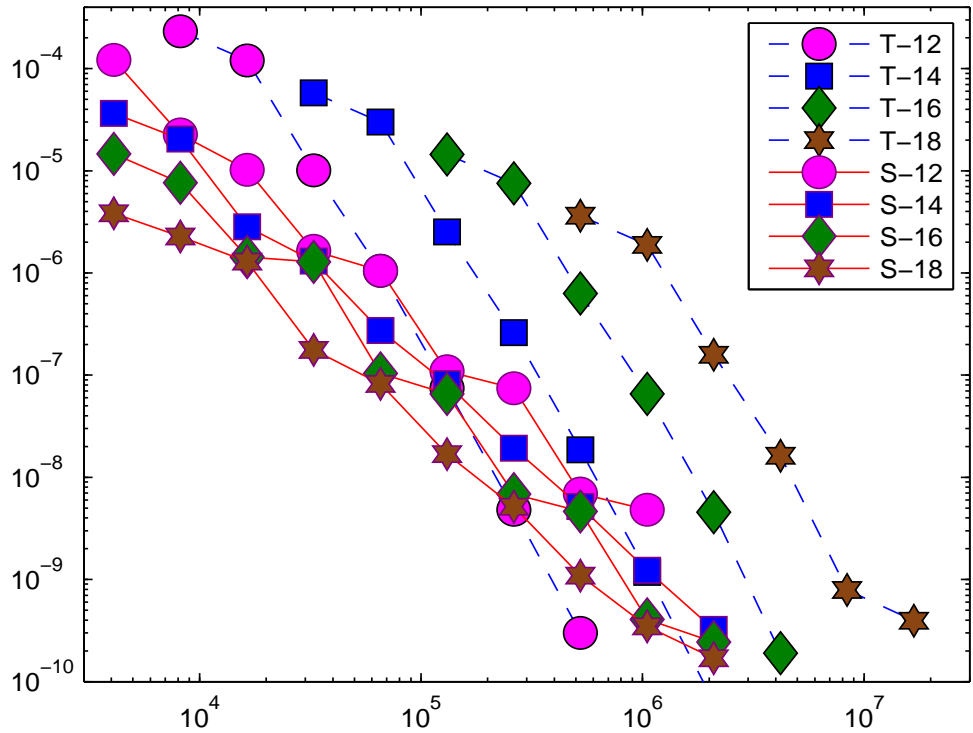,height=11cm,width=15cm}}
\caption{(NLS on the torus) Efficiency diagrams ({\em error} versus $N_{step}$) of SAM and TSFP4 for $\varepsilon=2^{-12},2^{-14},2^{-16},2^{-18}$.} 
\label{fig:efficient:per}
\end{figure}

Thirdly, we present efficiency diagrams. In Figure \ref{fig:efficient:per}, we represent error versus the total number $N_{step}$ of micro TSFP4 steps, for SAM and for TSFP4. Here, the red curves with the 'S' label plot errors of SAM and the dashed blue curves with the 'T' label plot errors of TSFP4. The error of SAM for a fixed number of TSFP4 steps is chosen as the minimal among all possible choices ($H_j,h_k$), i.e. $H_j=\pi/4/(2^j\eps), h_k = \pi/4/2^k$ with $j+k$ being fixed. Let us estimate from \eqref{error2} and \eqref{errTSFP4} theoretical values of the error versus $N_{step}$, for both methods. For SAM, we have $N_{step}=\frac{C_0}{hH\eps}$ so that, after a simple optimization on \eqref{error2}, one obtains an optimal error (assuming that the best choice for $H$, $h$ is ensured) under the form 
$$\mbox{\em error}_{\mbox{\scriptsize SAM}}\approx C_1\eps^4+\frac{C_0^2\sqrt{C_2C_3\eps}}{(N_{step})^2}.$$
For TSFP4, the number of micro-steps is $N_{step}=\frac{C_5}{\eps h}$. As a consequence, one gets from \eqref{errTSFP4}
$$\mbox{\em error}_{\mbox{\scriptsize TSFP4}}\approx \frac{C_4C_5^4}{\eps^3(N_{step})^4}.$$
On the curves in Figure \ref{fig:efficient:per}, we observe the expected behaviors. The error for SAM is proportional to $1/(N_{step})^2$ and improves proportionally to $\sqrt{\eps}$ when $\eps$ decreases, whereas the error for the TSFP4 method is proportional to $1/(N_{step})^4$ but is degraded as $\eps^{-3}$ when $\eps$ decreases. We can conclude that  TSFP4 performs better than SAM when $\varepsilon$ is relatively large,  e.g. $\varepsilon>2^{-12}$, and that SAM tends to perform better for smaller values of $\varepsilon$. The most striking result concerns the case $\eps=2^{-18}$ where the speed-up factor is 100.

\subsection{NLS on the torus: geometric properties} \label{sect:torusgeo}
In this subsection, we focus on geometric properties of SAM and thus consider  NLS on the torus \eqref{micro:periodic}, \eqref{micro:periodic:ini} with much longer time-intervals, of the form  $0\leq t\leq T_0/\eps^2$. Here, we are not so much interested by the accuracy of the wavefunction $\psi^\eps$ itself, but rather by the preservation of mass and energy, two invariants of \eqref{micro:periodic} respectively defined by
$$m^\eps(t)=\int_0^{2\pi}|\psi^\eps(t,x)|^2dx$$
and
$$\mathcal E^\eps(t)=\frac{1}{2}\int_0^{2\pi}|\partial_x\psi^\eps(t,x)|^2dx+\frac{\eps}{4} \int_0^{2\pi}2\cos(2x)|\psi^\eps(t,x)|^4dx.$$
The number of points for space discretization is taken as $N_x=32$. For the micro-integrator, we use the second-order Strang splitting scheme (TSFP2) with a micro time-step $h=2\pi/512$, ensuring that the CFL condition $h(N/2)^2<2\pi$ is satisfied so as  to avoid the effect of resonances over long times. Finally, we use the second-order formula \eqref{eq:Fh} for the interpolation of the vector field.

\begin{figure}[!htbp]
\centerline{\psfig{figure=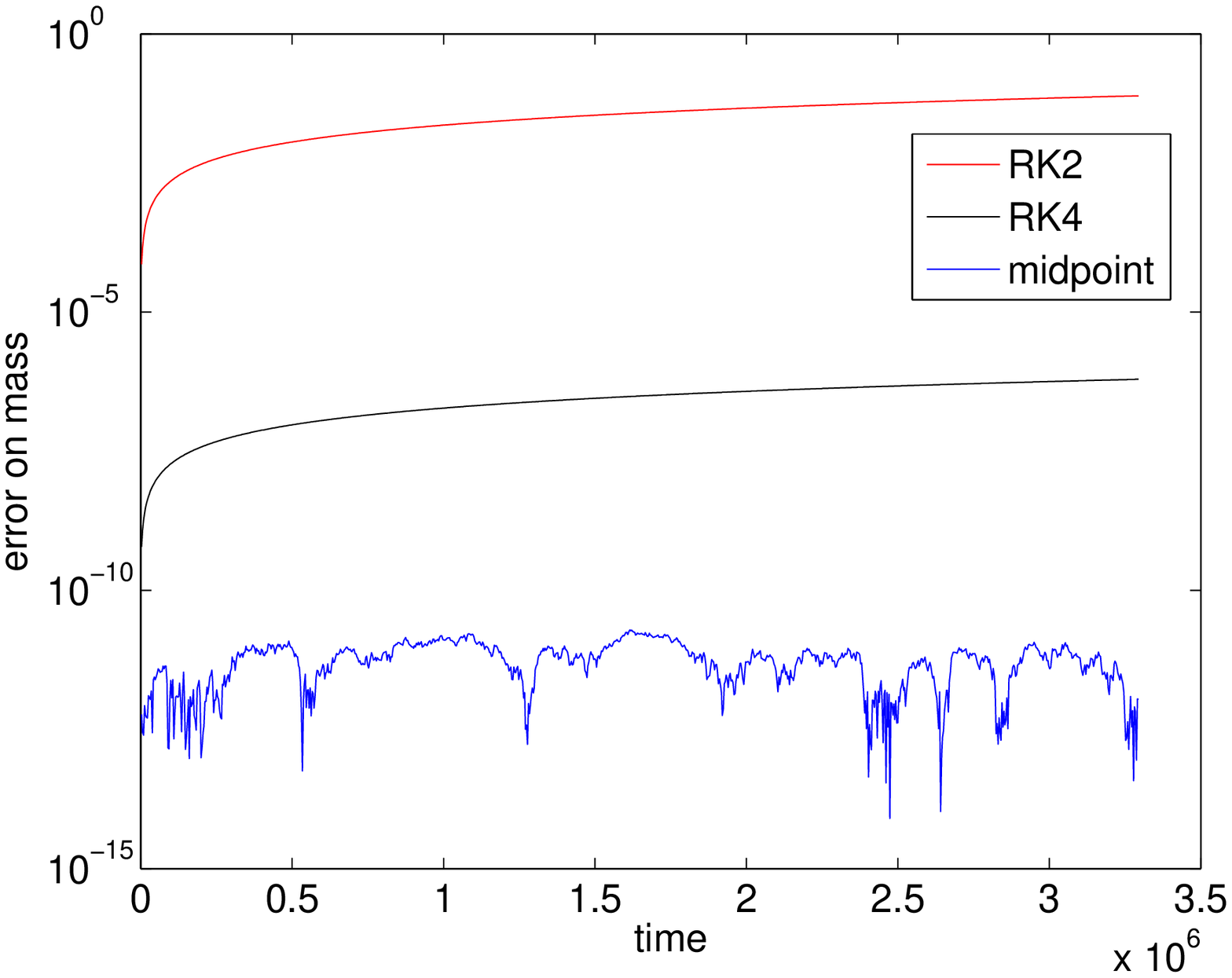,height=7cm,width=7.5cm}\hspace*{-3mm}
\psfig{figure=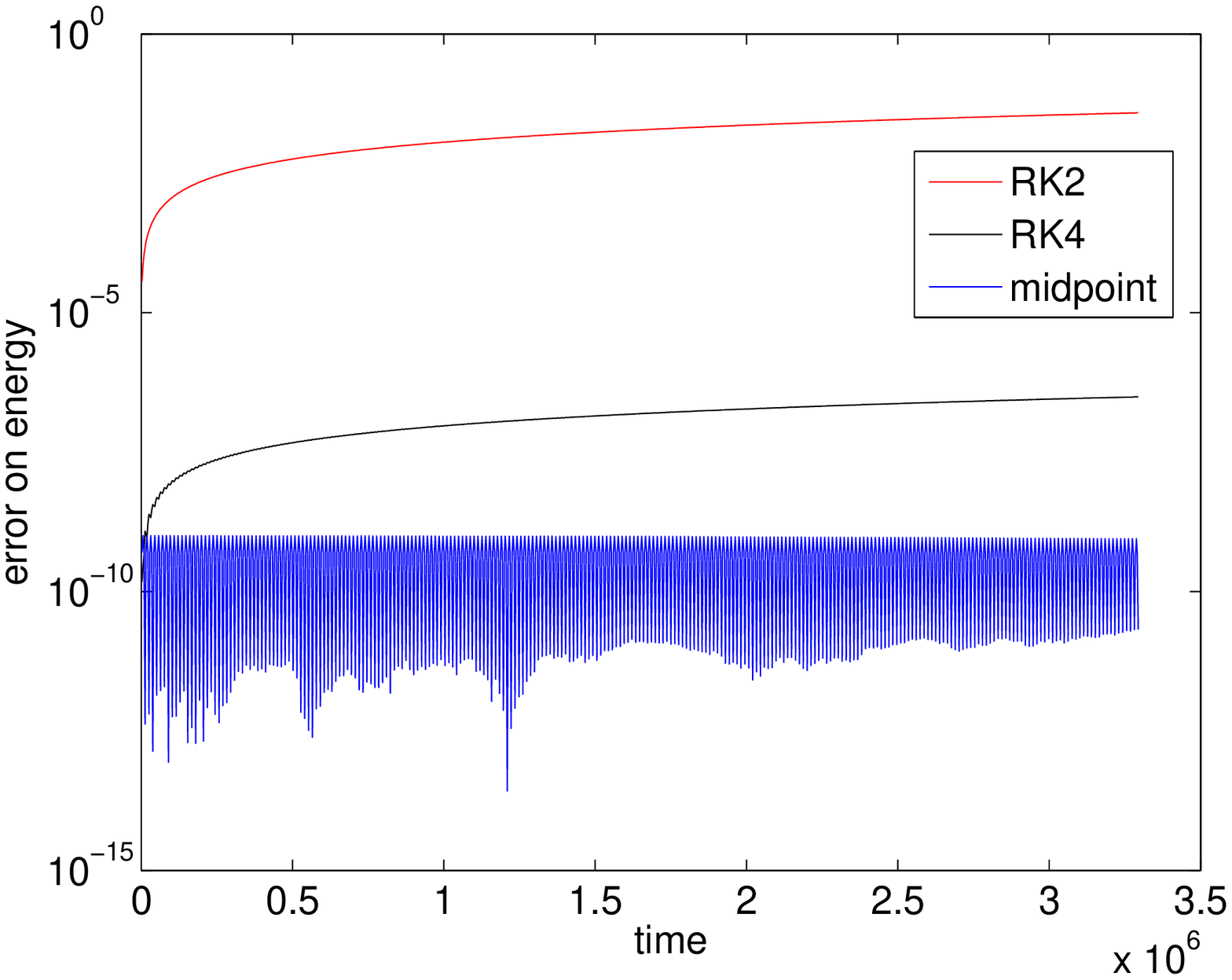,height=7cm,width=7.5cm}}
\caption{(NLS on the torus) Long time computation of two invariants with three macro-integrators: RK2, RK4 and the implicit midpoint scheme. Error (in log-scale) on mass (left) and error on energy (right) as a function of time $0\leq t\leq T_0/\eps^2$.} 
\label{fig:invariant}
\end{figure}

As for the macro-integrators, we experiment with three explicit Runge-Kutta methods, namely RK2, RK4 and the implicit midpoint rule, whose Butcher tableaux have the following values:
$$ \begin{array}{c|cc}
0\\
1/2 & 1/2\\
\hline
& 0 & 1
\end{array}, \qquad\qquad 
\begin{array}{c|cccc}
0\\
1/2 & 1/2\\
1/2 & 0 & 1/2\\
1 & 0 & 0 &  1\\
\hline
& 1/6 & 2/6 & 2/6 & 1/6
\end{array}, \quad \quad 
\begin{array}{c|c}
1/2 & 1/2 \\
\hline
& 1
\end{array}. $$ 
Note that the implicit midpoint rule is a symmetric method of order 2.
In our simulations, $\eps$ is given the value $2^{-11}\approx 5\times10^{-4}$ and the final integration time is now taken as $T_0/\eps^2$ with $T_0=\pi/4$. 
On Figure \ref{fig:invariant}, we plot in logarithmic scale the time evolution of the errors on mass (left) and energy (right), for the three macro-integrators considered. The macro time-step is $H=\pi/(2^7\eps)$. For the midpoint scheme, as expected, we observe a very good conservation of mass (error smaller than $2\times10^{-11}$) and energy (error smaller than $10^{-9}$), with no drift over this very long time interval. In contrast, the two non-symmetric schemes RK2 and RK4 display a linear drift in time: at the end of the simulation the errors on mass and energy are of order $10^{-2}$ for RK2 and of order $10^{-7}$ for RK4. These numerical results corroborate the statements of Subsection \ref{sect:geonum}.

%
%
%
\subsection{Gross-Pitaevskii equation}
\label{sect:GP}
In this subsection, we present our numerical results in the case (ii), namely NLS with a harmonic oscillator on the whole space,
\begin{align}\label{hermite_1d_micro_sch}
i\;  \partial_t\psi^\varepsilon&=  \left[-\frac{1}{2} \partial_{xx} + \frac{1}{2}( x^2-1)\right]  \psi^\varepsilon +\varepsilon  |\psi^\varepsilon|^2 \psi^\varepsilon,\qquad 0 \le t\le T_0/ \varepsilon,\quad x\in \R,\\
\label{hermite_1d_micro_sch_ini}\psi^\varepsilon(0,x) &= \psi_0(x), \quad x\in{\mathbb R}.
\end{align}
The micro-integrator is the fourth-order time-splitting Hermite spectral method (TSHP4)~\cite{BaoShen,Thalhammer}. The wavefunction is discretized in space by Hermite series as
$\psi(x) = \sum_{k=0}^N \widehat{\psi}_k\;h_k(x)$ where $h_k(x)$ ($k\in \N$), given explicitly in \cite{ShenJieBook3},  is the $(k+1)$-th eigenfunction of harmonic oscillator $L=-\frac{1}{2}\partial_{xx} +\frac{1}{2}( x^2-1)$ and is such that $Lh_k(x) = kh_k(x)$. In practice, we take $N=79$.  The macro-integrator for SAM is the RK4 scheme and the fourth-order interpolation formula is used for the approximation of the vector field.

The initial datum is $\psi_0(x)=h_0(x)+h_1(x)$ and the final integration time is $T_0/\eps$ with $T_0=2\pi$. For $\varepsilon = 2^{-3},\ldots,2^{-7}$, the reference solution $\psi^{ref}$ is obtained by the TSHP4 method with a time-step $h= \pi/10^3$ , while for $\varepsilon = 2^{-8},2^{-9},2^{-10}$ (and smaller $\varepsilon$), the reference solution is obtained  by SAM  with $H =\pi/(2^{10}\eps)$, $h = \pi/2^{10}$ and using 8th-order interpolation. As above, numerical errors are computed as discrete $\ell^2$ norm of the wavefunction
$$\mbox{\em error} =\|\psi^{ref}-\psi^n\|_{\ell^2} := \left(\sum_{j=0}^{N}|\psi^{ref}_j-\psi^n_j|^2 \omega_j\right)^{1/2},$$
where $\psi^n$ is the numerical solution and the coefficients $\omega_j$ are the rescaled Gauss-Hermite quadrature weights \cite{ShenJieBook3}.

\begin{figure}[!htbp]
\centerline{\psfig{figure=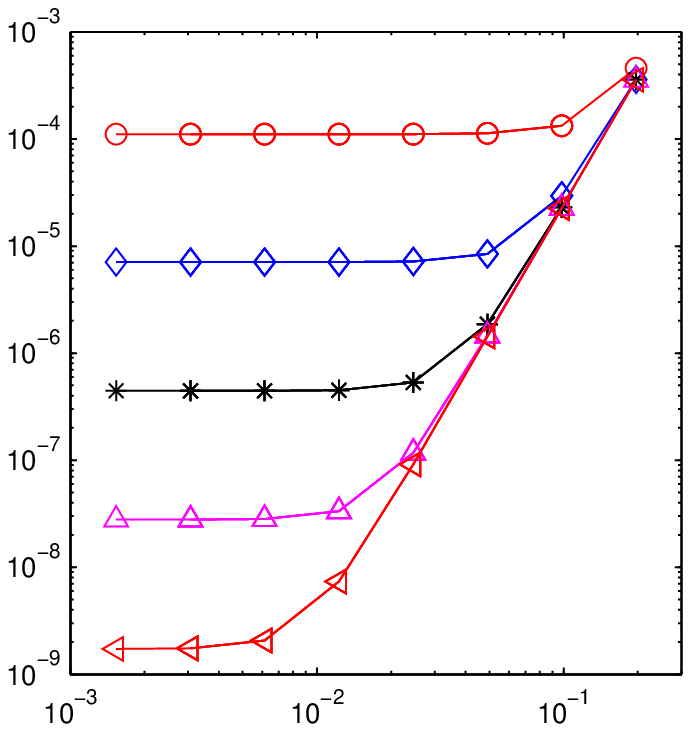,height=7cm,width=7.5cm}
\psfig{figure=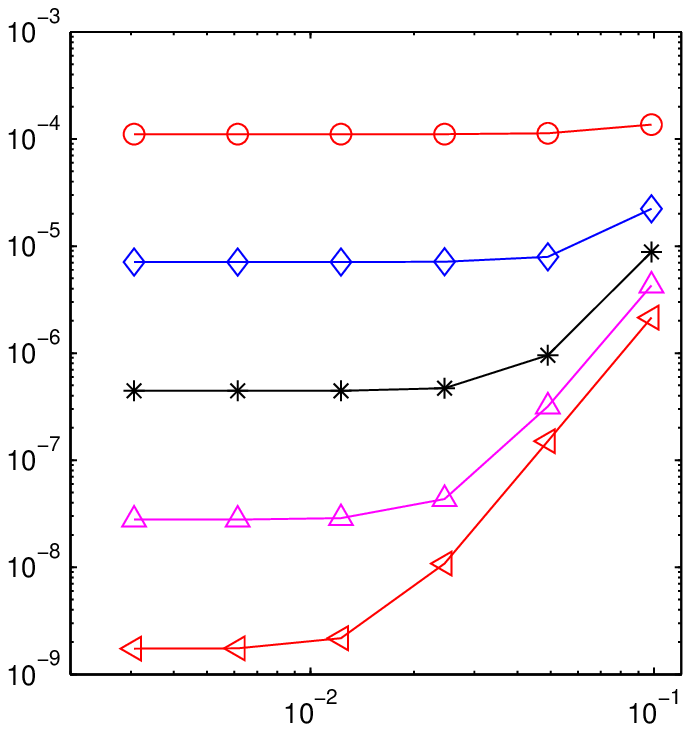,height=7cm,width=7.5cm}}
\caption{(Gross-Pitaevskii) Accuracy diagrams of SAM.  Error versus macro step $H$ (left) and error versus micro step $h$. For both diagrams, the values of  $\eps$ are $2^{-j}$, $j\in\{5,6,7,8,9\}$ (top to bottom).} 
\label{fig:rk4_tsfp4_hermite}
\end{figure}

In Figure \ref{fig:rk4_tsfp4_hermite}, we present accuracy diagrams of the error versus the macro step H (for different macro steps $H=\pi/(2^{j}\eps)$, $j = 3,4,\ldots,11$ with fixed micro steps $h=\pi/2^{12}$, and diagrams for the error versus the micro step $h =\pi/2^{5},\ldots,\pi/2^{9}$ with fixed macro steps $H=\pi/(2^{12}\eps)$.  The values $\varepsilon =2^{-5},\ldots,2^{-9}$ are tested. As in the case of NLS on the torus, these curves corroborate the error estimate \eqref{error2}.

In Table \ref{tab:tau:tshp}, we list the errors for the TSHP4 method as functions of $\eps$ and the micro step $h$. Again, our results also confirm the estimate
\begin{equation}
\label{errTSHP4}
\mbox{error}_{\mbox{\scriptsize TSHP4}}\approx C_4\varepsilon h^4.
\end{equation}

\begin{table}[t!]
\tabcolsep 0pt \caption{ Errors of TSHP4 for Gross-Pitaevskii.}
\label{tab:tau:tshp}\vspace*{-1em}
\begin{center}
\def\temptablewidth{1\textwidth}
{\rule{\temptablewidth}{1pt}}
\begin{tabular*}{\temptablewidth}{@{\extracolsep{\fill}}lllllll}
&$\varepsilon=2^{-3}$&$\varepsilon=2^{-4}$&$\varepsilon=2^{-5}$&$\varepsilon=2^{-6}$&$\varepsilon=2^{-7}$\\
$h=\pi/2^5$&1.16E-04&6.29E-05&3.28E-05&1.68E-05&8.47E-06\\
$h=\pi/2^6$&8.06E-06&4.37E-06&2.28E-06&1.17E-06&5.89E-07\\
$h=\pi/2^7$&5.18E-07&2.80E-07&1.46E-07&7.46E-08&3.73E-08\\
$h=\pi/2^8$&3.24E-08&1.75E-08&9.06E-09&4.45E-09&1.94E-09\\
\end{tabular*}
{\rule{\temptablewidth}{1pt}}
\end{center}
\end{table}


\section{Applications} \label{sect:appli}
In this section, we illustrate how SAM can be used to explore the qualitative behavior of highly-oscillatory systems. In the context of nonlinear Schr\"odinger equations of the form \eqref{NLS}, the dynamics of the coefficients (or modes) of the solution --when expanded on the basis of eigenfunctions of $A$-- exhibit interesting phenomena: The nonlinearity in \eqref{NLS} indeed induces coupling effects and energy transfers between modes, that we wish to visualize. As a first example, the 1D NLS equation on the torus studied in Subsection \ref{sect:torus} was already simulated in \cite{SAM} by SAM and other averaging techniques. Here, we simulate three other models: Gross-Pitaevskii equation in one dimension  and two models in two dimensions derived in the context of anisotropic confinement.

\bigskip
Let us briefly describe the phenomenology. The essential assumption here is that the unbounded operator $A$ has a compact resolvent and that the linear group $e^{-itA}$ is periodic in time. We recall that, if $\psi^\eps$ solves \eqref{NLS}, then $u^\eps=e^{itA}\psi^\eps$ satisfies the filtered equation
\begin{equation}
\label{eqfiltree2}
i\;\dot u^\eps(t) = \eps e^{it A} \left(\alpha(x)\left|e^{-itA} u^\eps \right|^2e^{-itA} u^\eps\right),\qquad u^\eps(0)=\psi_0.
\end{equation}
This equation is of the form \eqref{eq:HOP} and, up to exponentially small remainder terms, we have $u^\eps(t)=\Phi^\eps(t,\Psi^\eps(t,\psi_0))$, see \cite{SAM}. Recall that the change of variable $u\mapsto \Phi^\eps(t,u)$ is $P$-periodic in time and satisfies $\Psi^\eps(kP,u)=u$ for all $k\in\N$. Over the typical long time $T/\eps$ of a simulation, it accounts for the high oscillations that one may wish to avoid. In this situation,  it is more appropriate to represent only the function $\Psi^\eps(t,\psi_0)$ which, as solution to an autonomous equation \eqref{eq:aver}, is slowly varying and gives access to the long time nonlinear dynamics. SAM was precisely designed with this objective.

After projection on the eigenmodes of $A$ (in practice, Fourier or Hermite modes), \eqref{NLS} takes the equivalent form of an infinite system of coupled ODEs. Denote respectively by $(\lambda_k)_{k\in I}$,  $(e_k(x))_{k\in I}$, where $I$ is the set of indices ($\Z$ or $\N$ in applications), the eigenvalues and eigenfunctions of $A$. If 
$$\psi^\eps(t,x)=\sum_{k\in I}\psi_k(t)\,e_k(x),$$
solves \eqref{NLS}, then each mode $\psi_k$ satisfies the equation
\begin{eqnarray}
i\frac{d\psi_k}{dt}=\lambda_k\psi_k+\eps\sum_{\ell,m,n}\gamma_{k,\ell,m,n}\,\psi_\ell\,\overline{\psi_m}\,\psi_n,\end{eqnarray}
with $\displaystyle{\gamma_{k,\ell,m,n}=\int_\Omega \alpha(x)\overline{e_k(x)}e_\ell(x)\overline{e_m(x)}e_n(x)dx}$.

 Introducing the unknown $\xi_k(t)=e^{it\lambda_k}\psi_k(t)$, we obtain the filtered system
\begin{equation}
\label{xi}
i\frac{d\xi_k}{dt}=\eps\sum_{\ell,m,n}\gamma_{k,\ell,m,n}\,e^{it(\lambda_k-\lambda_\ell+\lambda_m-\lambda_n)}\,\xi_\ell\,\overline{\xi_m}\,\xi_n\,.
\end{equation}

\bigskip
In the sequel, we plot and compare the evolution of the moduli $|\psi_k(t)|=|\xi_k(t)|$ of the modes obtained by simulating three models:
\begin{itemize}
\item (TSFP4/TSHP4) The original equation \eqref{NLS}, discretized by a time splitting method used over the whole time interval $[0,T_0/\eps]$. 
\item (SAM) The averaged equation \eqref{eq:aver} associated to \eqref{xi}, discretized by SAM. Here, the macro solver is RK4, the micro solver is TSFP4/TSHP4 and the fourth-order interpolation formula is used for the approximation of the vector field. We recall that --up to numerical errors-- SAM and TSFP4/TSHP4 will coincide at the stroboscopic times $t=kP$, $k\in \N$.
\item (FAM) The first-order averaged model obtained by replacing $F^\eps$ in \eqref{eq:aver} by its first term in its expansion $F^\eps(u)=F_1(u)+\eps F_2(u)+\ldots$, see Remark \ref{remark}.

This model reads
\begin{equation}
\label{FAM}
i\;\dot u^\eps(t) = \eps \frac{1}{P}\int_0^P e^{i \tau A} \left(\alpha(x)\left|e^{-i\tau A} u^\eps (t,x)\right|^2e^{-i\tau A} u^\eps(t,x)\right)d\tau
\end{equation}
and the corresponding system on the modes $\xi_k$ reads
\begin{equation}\label{xiave}
i\frac{d\xi_k}{dt}= ~\eps \sum_{\ell,m,n \in \Lambda_k}\!
\gamma_{k,\ell,m,n}\,\xi_\ell\,\overline{\xi_m}\,\xi_n\end{equation}
with $\displaystyle{\Lambda_k=\{(\ell,m,n):\lambda_k\!-\!\lambda_\ell+\lambda_m\!-\!\lambda_n=0\}}$.

\end{itemize}
The numerical method for the approximation of the FAM is in fact directly constructed from \eqref{FAM}: this ODE is discretized by the RK4 scheme, where the integral defining the vector field is computed by the rectangle quadrature formula (which is spectrally convergent since the integrand is $P$-periodic). For each integral discretization, we use 64 points.

\subsection{Gross-Pitaevskii in dimension one} 
\begin{figure}[!htbp]
\centerline{
\psfig{figure=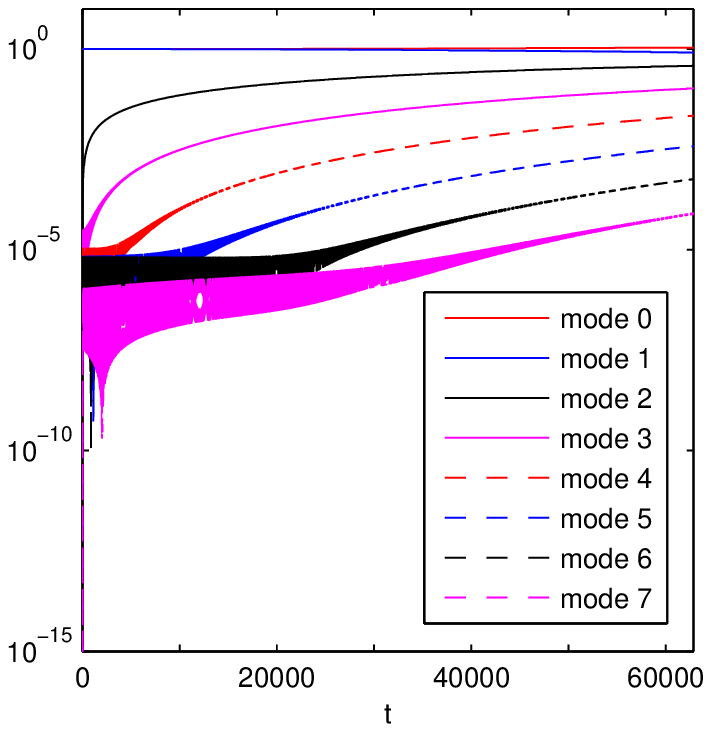,height=7cm,width=7.5cm}
 \psfig{figure=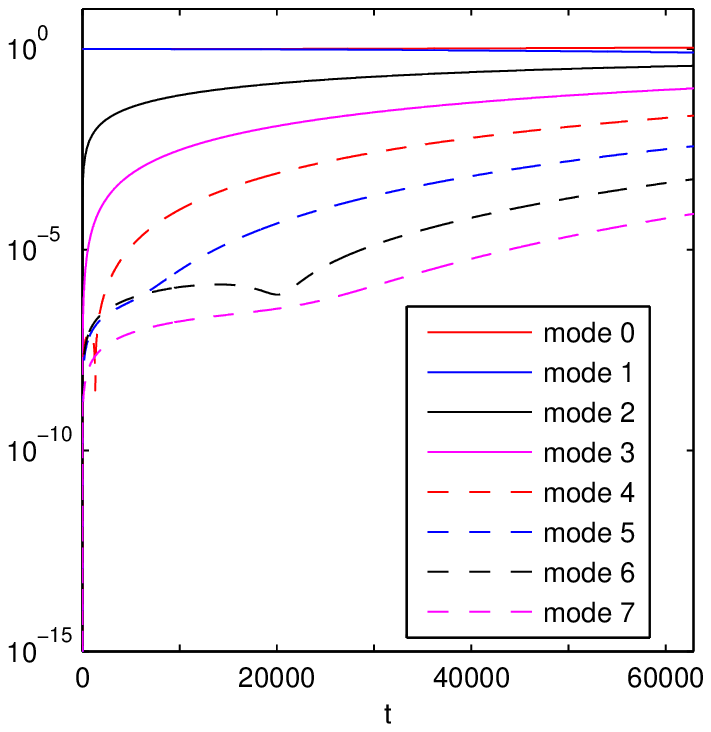,height=7cm,width=7.5cm}}
\centerline{\psfig{figure=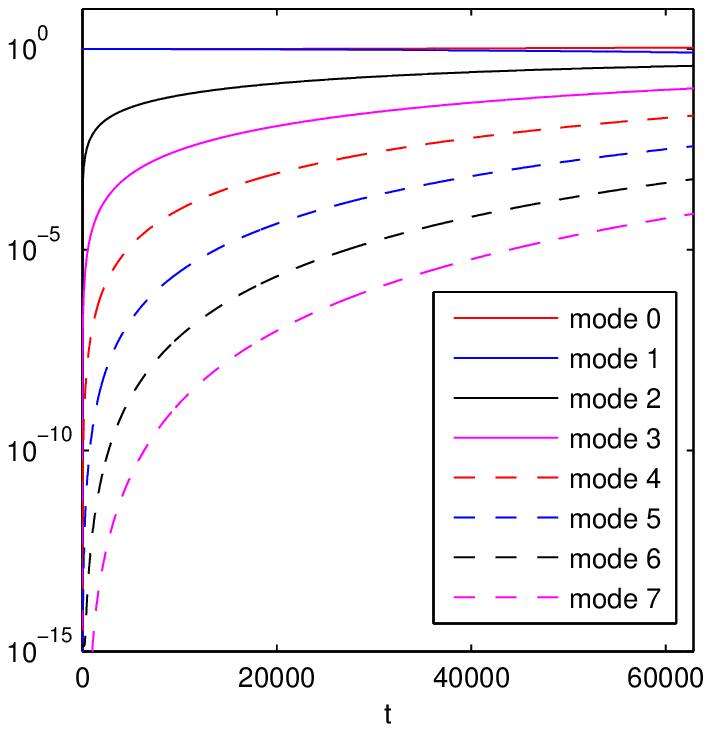,height=7cm,width=7.5cm}}
\caption{(Gross-Pitaevskii in dimension one) Evolution of Hermite coefficients by TSHP4 (left top), SAM (right top) and FAM (center bottom).} 
\label{fig:1d:herm}
\end{figure}

We consider here the one-dimensional equation \eqref{hermite_1d_micro_sch}, \eqref{hermite_1d_micro_sch_ini} with  $\varepsilon=10^{-4}$. The wave function is approximated by Hermite pseudospectral series with $N_x= 81$. Figure \ref{fig:1d:herm} plots the evolution of the absolute value of the Hermite coefficients by TSHP4, SAM and FAM. Numerical parameters are taken as follows: The TSHP4 method is applied with $h=\frac{2\pi}{10^3}$, SAM with $(H,h) = (\frac{2\pi}{10^4\eps },\frac{2\pi}{400})$ and FAM with $h = \frac{2\pi}{10^4}$. 

We observe an interesting energy cascade phenomenon, very similar to the case of NLS on the torus $[0,2\pi]\times [0,2\pi]$ studied in \cite{cf} (see also \cite{SAM} for numerical simulations by SAM). While the modes greater than 1 are equal to zero initially, they grow and become significantly large in a characteristic time that depends on the mode. The point is that, as long as a mode is of order $\mathcal O(\eps)$, it is highly-oscillatory and its observation by TSFP4 is not convenient on the first diagram (Fig. \ref{fig:1d:herm} top left). Instead, as we said in the introduction of this section, SAM cleans up the oscillations associated to the change of variable $\Phi^\eps$ and it is finally easier to observe the energy cascade on the second diagram (Fig. \ref{fig:1d:herm} top right). The curves represented on third diagram (FAM) are also very smooth, but show a dynamics that is correct only when modes have reached a value above a threshold $\mathcal O(\eps)$: under this level, the dynamics given by FAM is inaccurate.

\subsection{2D NLS on an anisotropic torus}
We consider the cubic NLS equation on a very thin bidimensional torus:
\begin{align*}
i\partial_t \psi^\varepsilon &= -\Delta \psi^\eps+ |\psi^\varepsilon|^2 \psi^\varepsilon, \qquad 0 < t \le T_0,\quad (x,y)\in [0,2\pi]\times [0,2\pi \eps],\\
\psi^\varepsilon(0,x,y) &= \psi_0(x,y/\eps),\qquad (x,y)\in [0,2\pi]\times [0,2\pi \eps].
\end{align*} 
Here, the dimensionless parameter $\eps$ is the ratio between the typical length-scales in directions $y$ and $x$. Let us rescale the variables, setting $(t,x,y)=(\eps^2t',x',\eps y')$ (then omitting the primes for simplicity in the new equation). The space domain $\T_{2\pi}^2=[0,2\pi]\times [0,2\pi]$ is now independent of $\eps$ and we obtain
\begin{align}\label{micro:eqn}
i\partial_t\psi^\varepsilon&=  - \partial_{yy}\psi^\eps +\varepsilon^2 \left(-\partial_{xx} +  |\psi^\eps|^2 \right)\psi^\eps,\qquad  0 < t \le T_0/ \varepsilon^2,\\
\psi^\varepsilon(0,x,y) &= \psi_0(x,y),\qquad (x,y)\in \T_{2\pi}^2.
\end{align}  
This equation is still of the form \eqref{eq:aut}. However, the nonlinearity $g$ is now unbounded in any Sobolev space. In such case, the exponential estimates proved in \cite{SAM} are not valid. One can nevertheless expect polynomial estimates for smooth initial data, which are enough to give foundations to SAM.

\begin{figure}[!htbp]
\centerline{
\psfig{figure=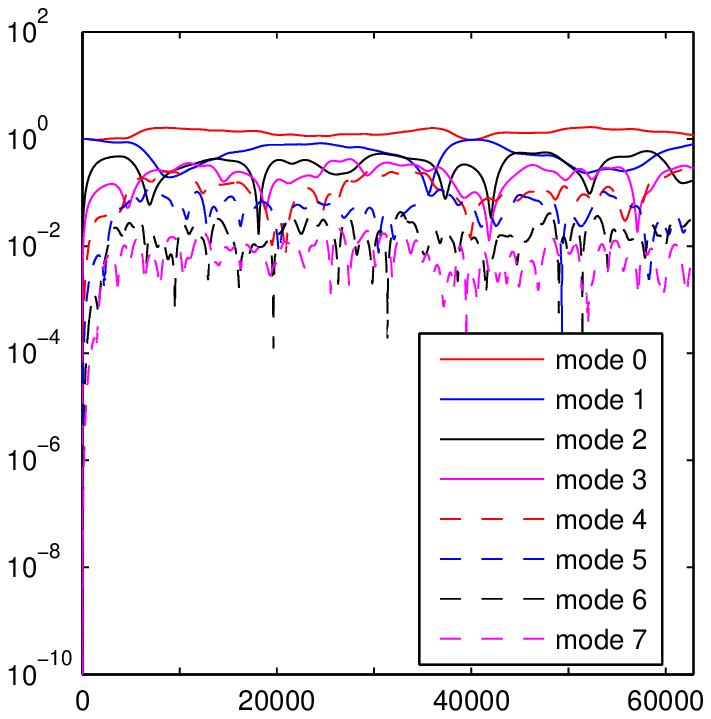,height=7cm,width=7.5cm}
\psfig{figure=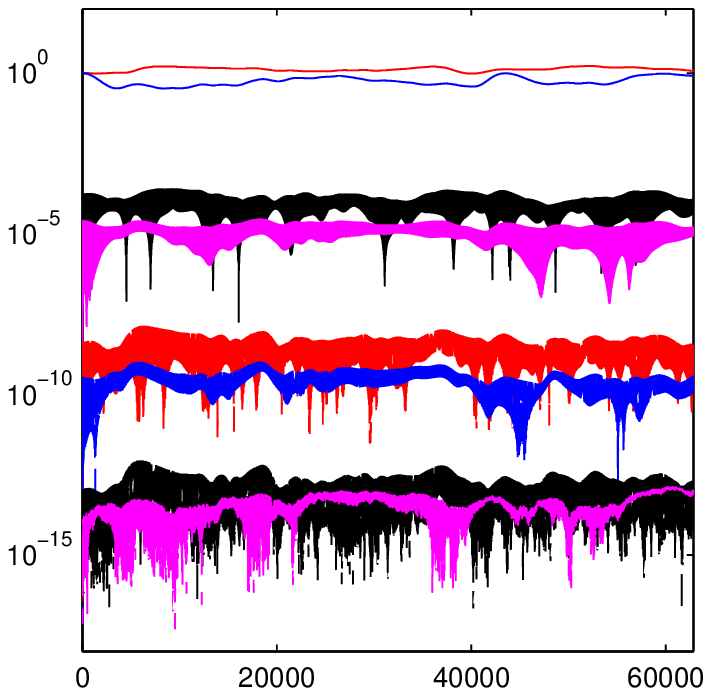,height=7cm,width=7.5cm}}

\centerline{
\psfig{figure=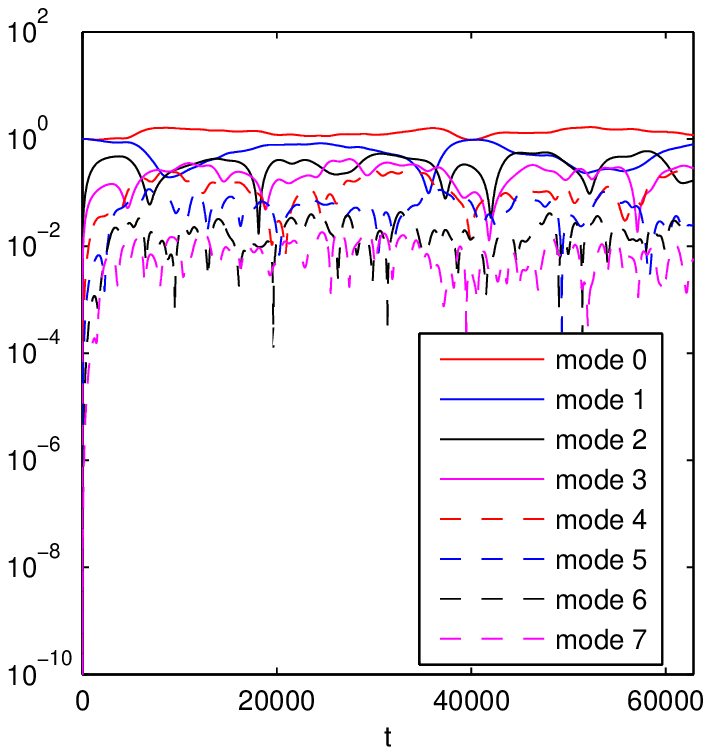,height=7cm,width=7.5cm}
\psfig{figure=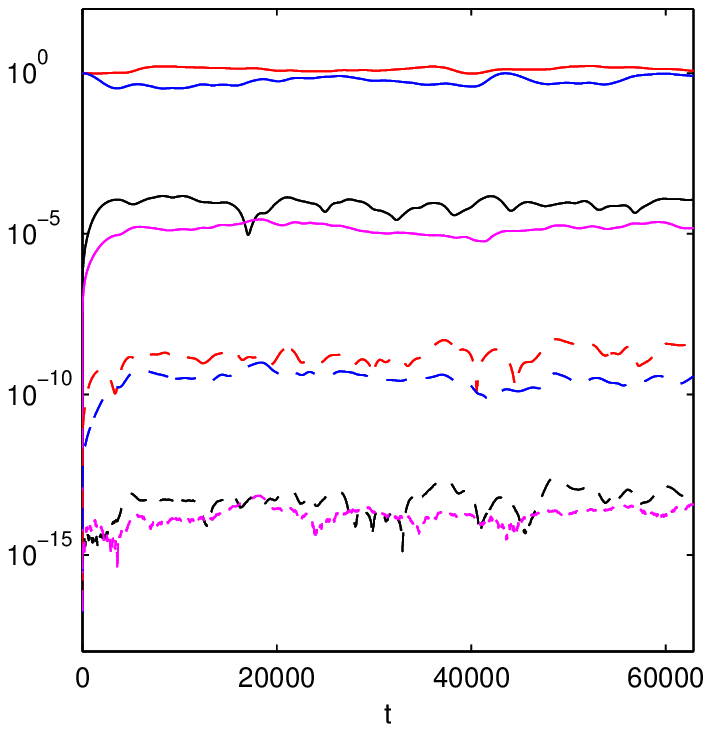,height=7cm,width=7.5cm}}

\centerline{
\psfig{figure=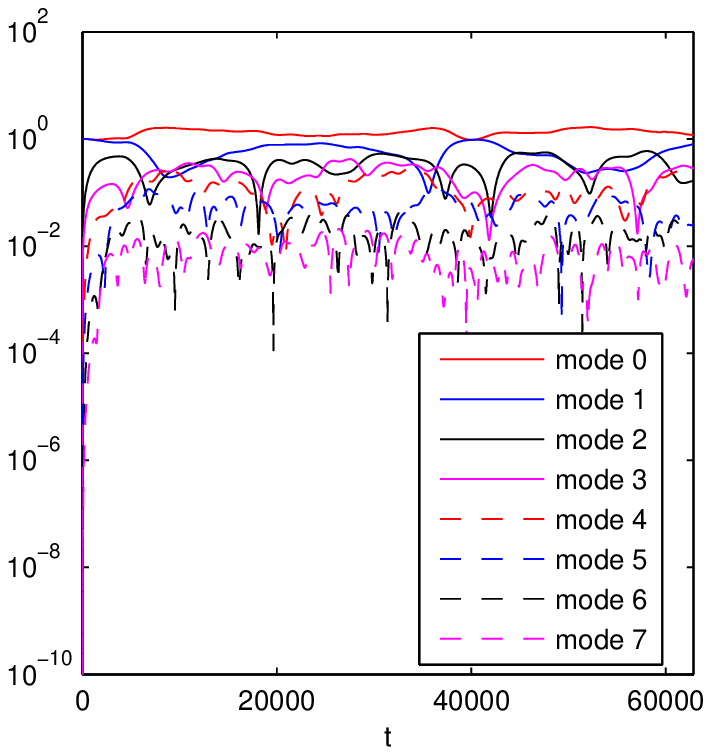,height=7cm,width=7.5cm}
\psfig{figure=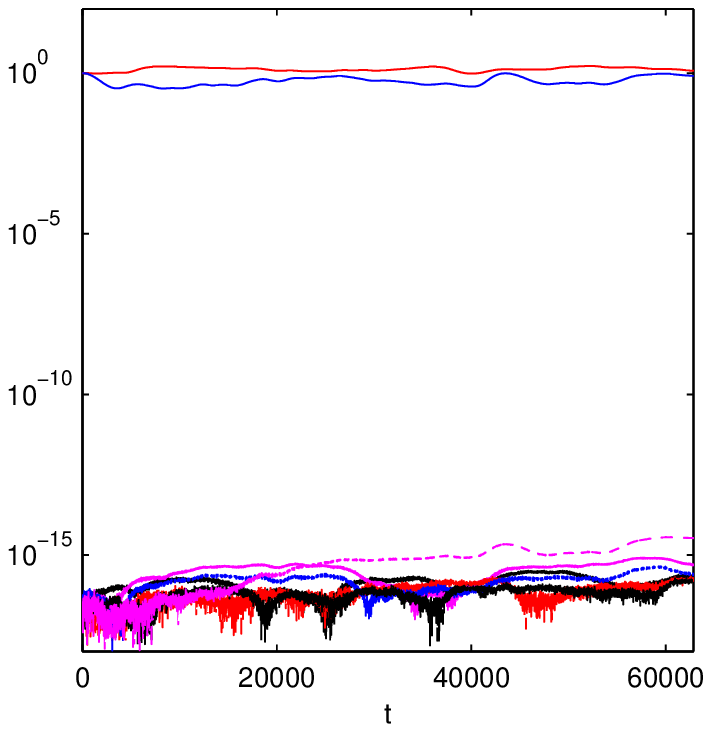,height=7cm,width=7.5cm}}

\caption{(2D NLS on an anisotropic torus) Evolution of Fourier mode by TSFP4(Top), SAM (middle) and FAM (bottom):  $|\widehat{\psi}_{n,0}|$ (left)  and $|\widehat{\psi}_{0,n}|$(right).} 
\label{period:2d:mode}
\end{figure}

As initial data for \eqref{micro:eqn}, we take $\psi_0(x,y)=1+2\cos(x)+2\cos(y)$.
The numerical parameters of our simulations are the following. The wave function is approximated by Fourier series with $N_x = 128$ and $N_y = 128$. The TSFP4 method is applied with a time step $h= \frac{2\pi}{500} $, 
SAM is applied with a macro step $ H = \frac{2\pi}{10^4\eps^2}$ and a micro step $h = \frac{2\pi}{400}$ and the first averaged model (FAM) is computed with the time step $h=\frac{2\pi}{10^4}$. In Fig. \ref{period:2d:mode}, we plot the evolution of the Fourier coefficients $|\widehat{\psi}_{n,0}|$  and $|\widehat{\psi}_{0,n}|$ with $n=0,1,\ldots,7$, obtained by the three different methods, for $\varepsilon=0.01$ and $t\in [0, 2\pi/\eps^2]$. 

The following observations can be made. The three diagrams on the left are very similar, indicating that the first averaged model FAM (which is cheaper) is sufficient to describe correctly the evolution of modes $|\widehat{\psi}_{n,0}|$. However, on the three diagrams of the right, it can be inferred  that FAM is not able to capture the dynamics of higher modes $k_y\geq 2$. Some interesting phenomena appearing on these higher modes in $y$ can only be observed by SAM ( on TSFP4, these modes are highly-oscillatory): for instance, the modes seem to arranged by pairs as
$$|\widehat{\psi}_{0,0}|,|\widehat{\psi}_{0,1}|=\mathcal O(1),\; |\widehat{\psi}_{0,2}|,|\widehat{\psi}_{0,3}|=\mathcal O(\eps^2),\; |\widehat{\psi}_{0,4}|,|\widehat{\psi}_{0,5}|=\mathcal O(\eps^4),\; |\widehat{\psi}_{0,6}|,|\widehat{\psi}_{0,7}|=\mathcal O(\eps^6),\,\ldots\quad \qquad \qquad $$

\subsection{2D anisotropic Gross-Pitaevskii equation}
We consider here the Gross-Pitaevskii equation with a strongly anisotropic harmonic confinement potential:
\begin{align*}
&i\partial_t \psi^\varepsilon = -\frac12\Delta \psi^\eps + \frac12\left(x^2+\frac{y^2}{\varepsilon^4}\right)\psi^\eps+ \beta | \psi^\varepsilon|^2 \psi^\varepsilon,\quad 0 < t \le T_0,\quad (x,y)\in \R^2,\\
&\psi^\varepsilon(0,x,y) = \psi_0(x,y/\eps),  \qquad  (x,y)\in {\mathbb R}^2,
\end{align*}
see \cite{CaiSeond} for the physical context.

In order to rewrite this equation in our framework, the adequate change of variable is again $(t,x,y)=(\eps^2t',x',\eps y')$, which yields the following model on the new unknown (still denoted $\psi^\eps(t,x,y)$):
\begin{align}
\label{GPconf}
&i\partial_t \psi^\varepsilon = \left(-\frac12\partial_{yy}+\frac{y^2}2\right)\psi^\eps+\eps^2 \left(-\frac12\partial_{xx}+\frac{x^2}{2}+\beta | \psi^\varepsilon|^2\right) \psi^\varepsilon,\; 0 < t \le T_0/\eps^2, \\
&\psi^\varepsilon(0,x,y) = \psi_0(x,y),  \quad  (x,y)\in {\mathbb R}^2.
\end{align}
Our initial data is
$$\psi_0(x,y)=h_0(y) (h_0(x)+h_2(x)),$$
where $h_k$ are the Hermite functions (see Subsection \ref{sect:GP}). The wave function is approximated by Hermite pseudospectral series with $N_x =N_y= 81$. In the computations, TSHP4 is applied with a time step $h = \frac{P}{100}$, SAM is used with $(H,h) =( \frac{P}{10^3\eps^2}, \frac{P}{100})$ and FAM is used with the time step $h = \frac{P}{10^3}$.

\begin{figure}[!htbp]
\centerline{
\psfig{figure=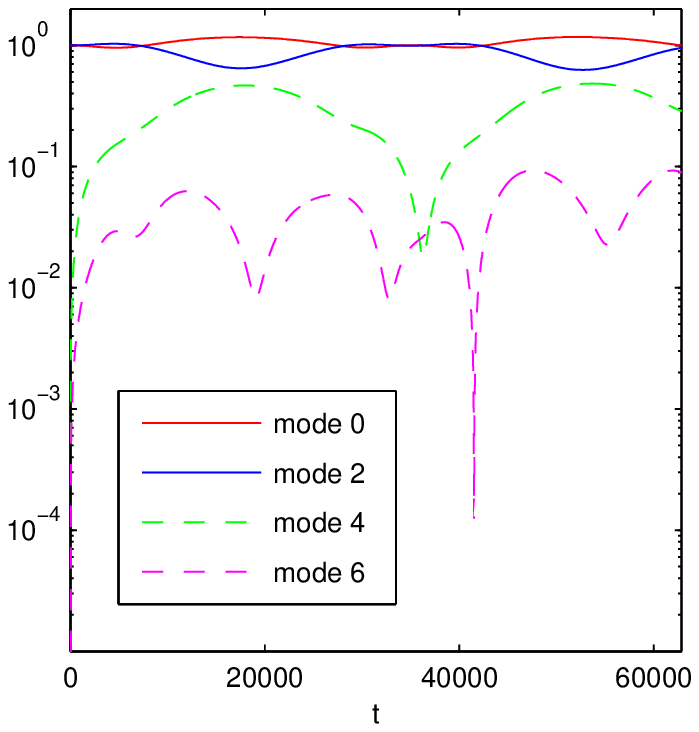,height=5cm,width=5.5cm}\hspace*{-5mm}
 \psfig{figure=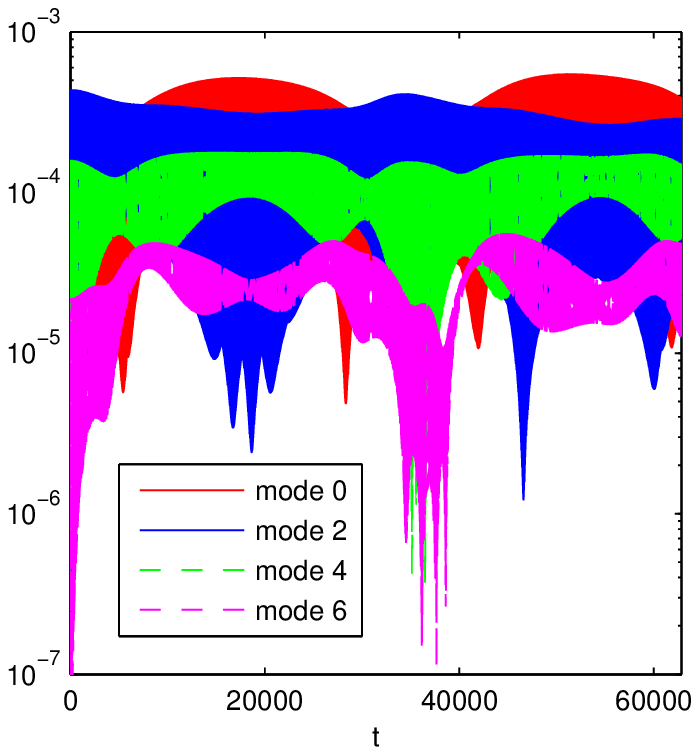,height=5cm,width=5.5cm}\hspace*{-5mm}
\psfig{figure=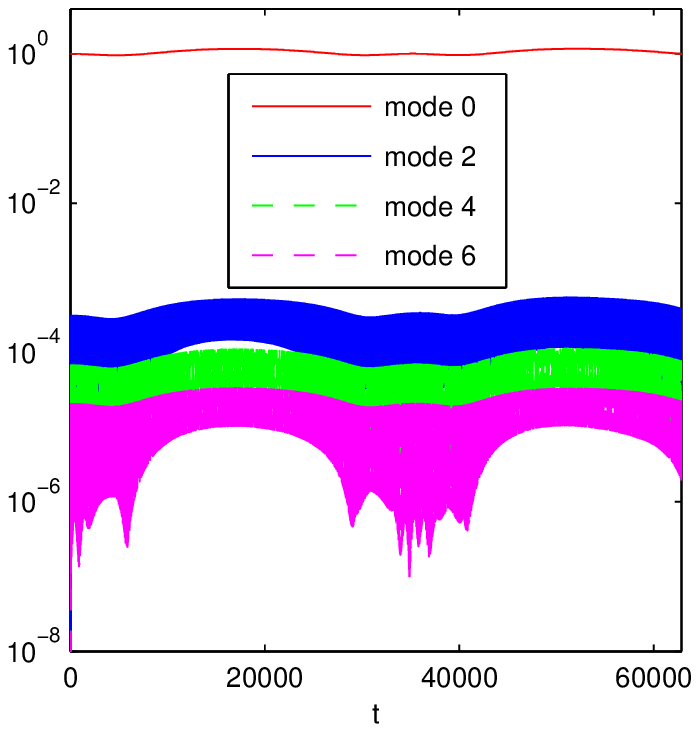,height=5cm,width=5.5cm} }

\centerline{
\psfig{figure=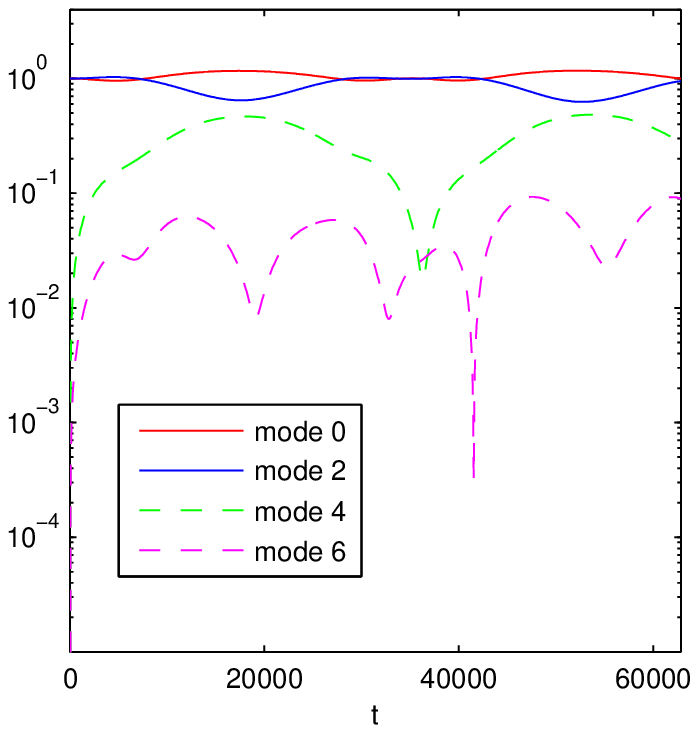,height=5cm,width=5.5cm}\hspace*{-5mm}
\psfig{figure=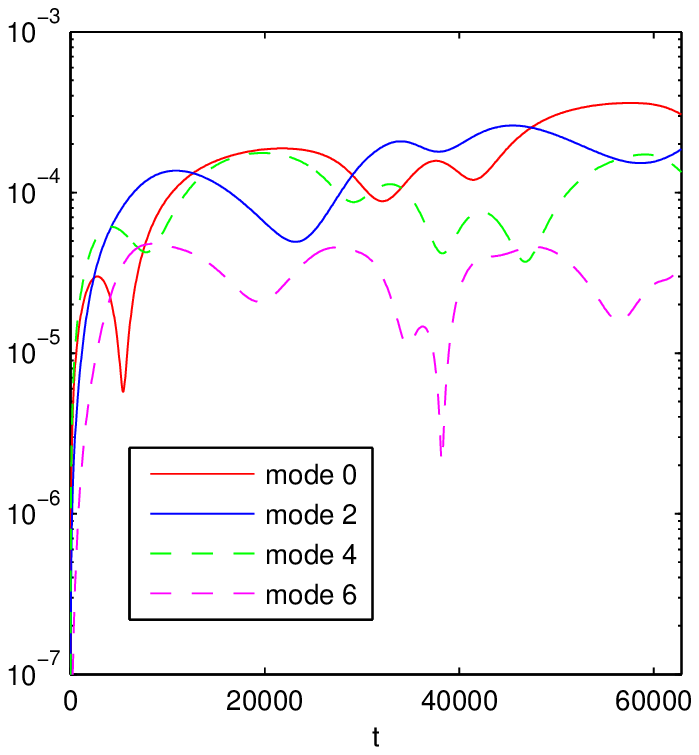,height=5cm,width=5.5cm}\hspace*{-5mm}
\psfig{figure=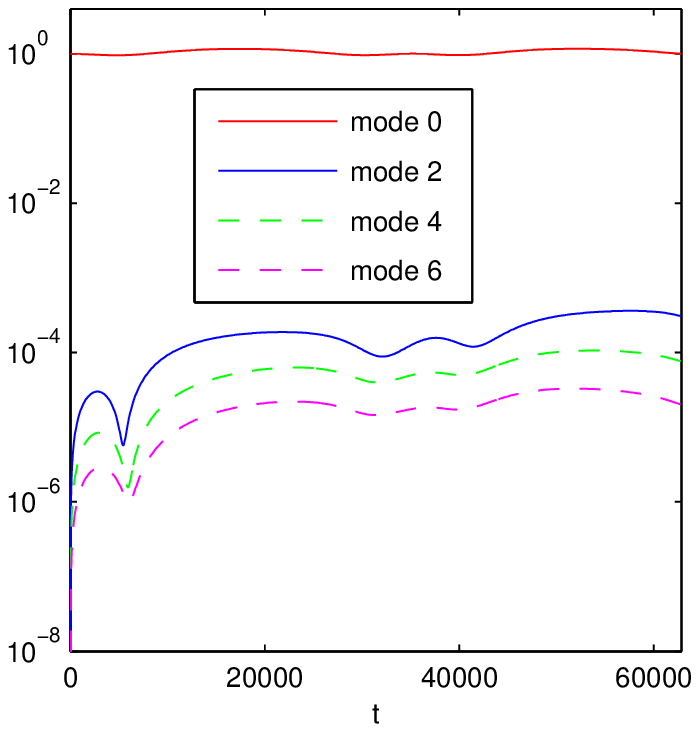,height=5cm,width=5.5cm}}

\centerline{ 
\psfig{figure=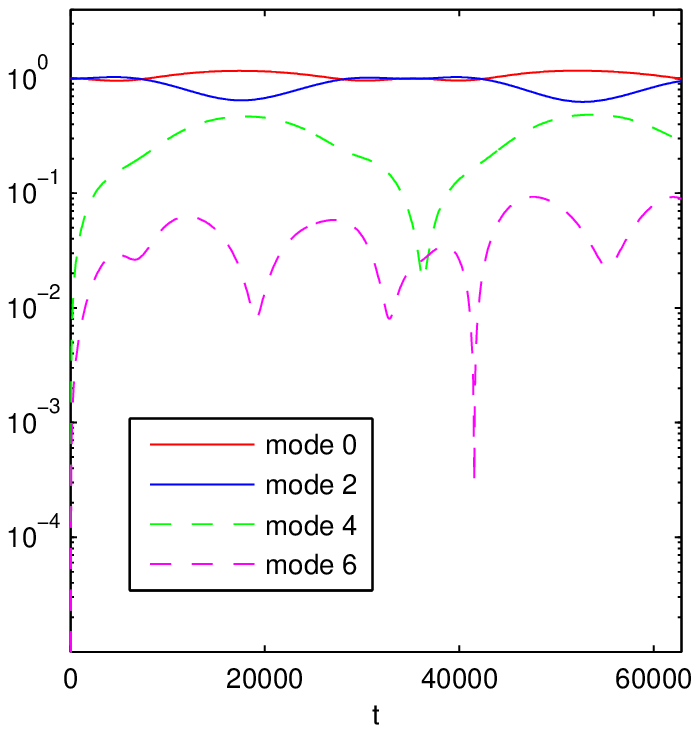,height=5cm,width=5.5cm}\hspace*{-5mm}
\psfig{figure=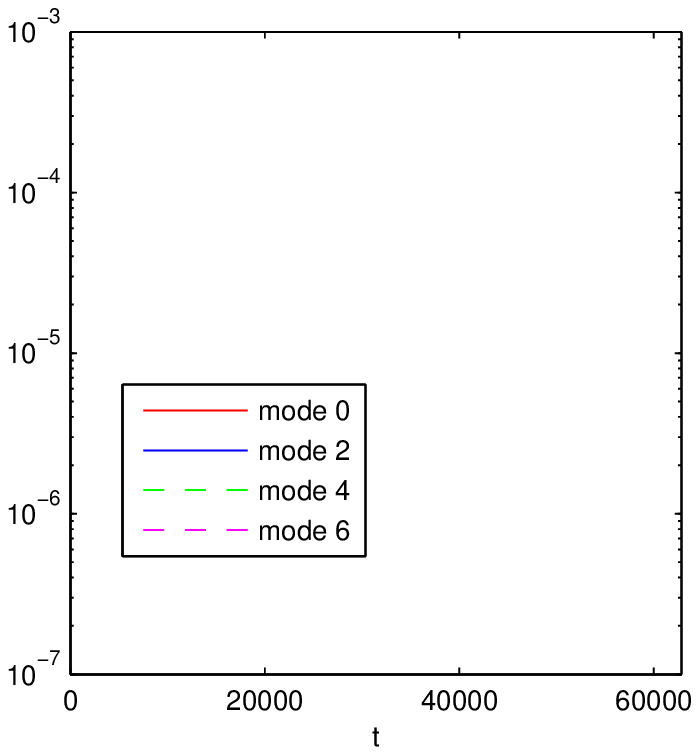,height=5cm,width=5.5cm}\hspace*{-5mm}
\psfig{figure=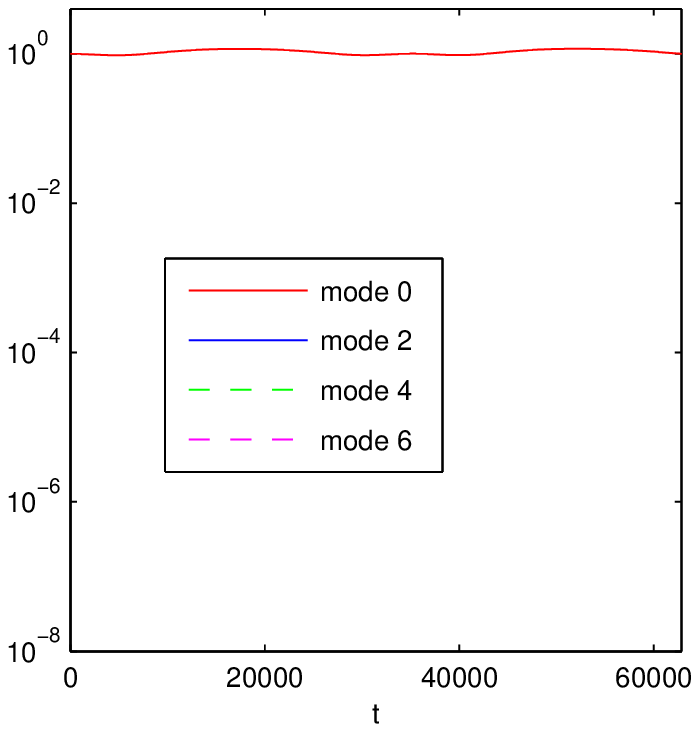,height=5cm,width=5.5cm}}

\caption{(2D anisotropic Gross-Pitaevskii equation) Evolution of Hermite modes by TSHP4(top), SAM (middle) and FAM (bottom):  $|\widehat{\psi}_{n,0}|$ (left), $|\widehat{\psi}_{n,2}|$ (middle) and $|\widehat{\psi}_{0,n}|$(right).} 
\label{hermite:2d:mode}
\end{figure}

On Figure \ref{hermite:2d:mode}, we plot the evolution of Hermite modes for $\varepsilon= 10^{-2}$ and $\beta=5$. Similar comments as for 2D NLS on an anisotropic torus can be done. For a description of modes 0 in $y$, FAM is sufficient. Interesting dynamics can be observed on higher modes in $y$, and SAM is the more appropriate model for the investigation of these dynamics, since it filters out all the oscillations.

\section*{Acknowledgement}
The computational results presented have been achieved using the Vienna Scientific Cluster (VSC).



\begin{thebibliography}{alpha}

\bibitem{BaoJM}
{\sc W.~Bao, S.~Jin and P.A.~Markowich}, {\em On Time-splitting spectral approximations for the \sch\ equation in the semiclassical regime}, J. Comput. Phys., 175 (2002), pp.~487--524.

\bibitem{BaoJM2}
{\sc W.~Bao, S.~Jin and P.A.~Markowich}, {\em Numerical studies of time-splitting spectral discretizations of nonlinear \sch\ equations in the semiclassical regime}, SIAM J. Sci. Comput., 25 (2003), pp.~27--64.

\bibitem{BaoShen}
{\sc W.~Bao  and J.~Shen}, {\em A fourth-order time-splitting Laguerre-Hermite pseudospectral method for Bose-Einstein condensates}, SIAM J. Sci. Comput., 26 (2005), pp.~2010--2028.

\bibitem{CaiSeond}
{\sc N. Ben Abdallah, Y.Y.~Cai, F.~Castella and F.~M\'ehats}, {\em Second order averaging for the nonlinear Schr\"{o}dinger equation with strongly anisotropic potential},
Kinetic and Related Models,  4 (2011), pp.~831--856.

\bibitem{TimeAveJDE}
{\sc N. Ben Abdallah, F.~Castella and F.~M\'ehats}, {\em Time averaging for the strongly confined nonlinear Schr\"{o}dinger equation, using almost-periodicity}, 
J. Differ. Equations, 245 (1) (2008), pp.~154--200.


\bibitem{BCOR09} 
{\sc S. Blanes, F. Casas, J.A. Oteo and J. Ros}, {\em The Magnus and expansion and some of its applications}, Physics Reports, 470 (2009), pp.~151--238.

\bibitem{CCMSS11b} 
{\sc M.P.~Calvo, P.~Chartier, A.~Murua and J.M.~Sanz-Serna}, {\em A stroboscopic numerical method for highly
oscillatory problems},  in Numerical Analysis and Multiscale Computations, Lect. Notes Comput. Sci. Eng., 82 (2012), pp.~73--87.

\bibitem{CCMSS11}
{\sc M.P.~Calvo, P.~Chartier, A.~Murua and J.-M.~Sanz-Serna}, {\em Numerical stroboscopic averaging for ODEs and DAEs},  Appl. Numer. Math., 61 (2011), pp.~1077--1095.

\bibitem{cf} 
{\sc R. Carles and E. Faou}, {\em Energy cascades for NLS on the torus}, Discrete and Continuous Dynamical Systems-Series A (DCDS-A), 32 (6)  (2012), pp.~2063--2077. 

\bibitem{SAM}
{\sc F.~Castella, P.~Chartier, F.~M\'ehats and A.~Murua}, {\em Stroboscopic averaging for the nonlinear Schr\"{o}dinger equations}, preprint.

\bibitem{superconvergence} 
{\sc P.~Chartier, F.~M\'ehats, M.~Thalhammer and Y.~Zhang}, {\em A note on the convergence of splitting methods for periodic highly-oscillatory systems}, preprint.

\bibitem{CMSS12b}
{\sc P.~Chartier, A.~Murua and J.M.~Sanz-Serna}, {\em A formal series approach to averaging: exponentially
small error estimates}, Discrete and Continuous Dynamical Systems-Series A (DCDS-A), 32 (9) (2012), pp.~3009--3027.

\bibitem{CMSS10} 
{\sc P.~Chartier, A.~Murua and J.M.~Sanz-Serna},{\em Higher-order averaging, formal series and
numerical integration I: B-series}, Found Comput Math, 10 (2010), pp.~695--727.

\bibitem{E03} 
{\sc W. E}, {\em Analysis of the heterogeneous multiscale method for ordinary differential
equations}, Comm. Math. Sci., 1 (2003), pp.~423--436.

\bibitem{E03b} 
{\sc W.~E and B.~Engquist}, {\em The heterogeneous multiscale methods}, Comm. Math. Sci., 1 (2003), pp.~87--132.

\bibitem{E07} {\sc W.~E, B.~Engquist, X.~Li, W. Ren and E. Vanden-Eijnden}, {\em Heterogeneous multiscale
methods: A review}, Commun. Comput. Phys., 2 (2007), pp.~367--450.

\bibitem{E05} {\sc B. Engquist and R. Tsai}, {\em Heterogeneous multiscale methods for stiff ordinary differential
equations}, Math. Comput., 74 (2005), pp.~1707--1742.

\bibitem{Grebert}
{\sc B.~Gr\'ebert and C.~Villegas-Blas}, {\em On the energy exchange between resonant modes in nonlinear Schr\"{o}dinger equations}, Ann. I. H. Poincar\'e, 28 (2011), pp.~127--134.

\bibitem{Hairer06} 
{\sc E. Hairer, C. Lubich and G. Wanner}, {\em Geometric Numerical Integration: Structure-Preserving
Algorithms for Ordinary Differential Equations}, Springer, 2006. 

\bibitem{neishtadt84}
{\sc A.I.~Neishtadt},{\em The separation of motions in systems with rapidly rotating phase}, J. Appl. Math. Mech., 48 (1984), pp.~133--139.

\bibitem{Sanders07}
{\sc J.A.~Sanders, F.~Verhulst and J.~Murdock}, {\em Averaging methods in nonlinear dynamical systems},
Springer, 2007. 

\bibitem{ShenJieBook3}
{\sc J.~Shen, T.~Tang and L.-L.~Wang}, {\em Spectral Methods: Algorithms, Analysis and Applications}, Springer, 2011.

\bibitem{Strang}
{\sc G. Strang}, {\em On the construction and comparision of difference schemes}, SIAM J. Numer. Anal., 5 (1968), pp.~505--517.

\bibitem{Thalhammer}
{\sc M.~Thalhammer, M.~Caliari and C.~Neuhauser}, {\em High-order time-splitting Hermite and Fourier spectral methods}, J.
Comput. Phys., 228 (2009), pp.~822-832.

\bibitem{Yoshida}
{\sc H. Yoshida}, {\em Construction of higher order symplectic integrators}, Phys. Lett. A, 150 (1990), pp.~262--268.

\bibitem{SPMCompare}
{\sc Y. ~Zhang and X.~Dong}, {\em On the computation of ground state and
dynamics of Schr\"{o}dinger-Poisson-Slater system}, J. Comput. Phys.,  230 (2011), pp.~2660--2676.

\end{thebibliography}
\end{document}